\documentclass[11pt]{amsart}
\usepackage{amsmath,amssymb}
\usepackage[matrix,arrow,curve,frame]{xy}    
                                                                               
\xymatrixcolsep{1.9pc}                          
\xymatrixrowsep{1.9pc}
\newdir{ >}{{}*!/-5pt/\dir{>}}                  

\usepackage{mathrsfs}
\usepackage{hyperref}

\swapnumbers

\newtheorem{thm}{Theorem}[section]

\newtheorem{prop}[thm]{Proposition}
\newtheorem*{clm}{Claim}
\newtheorem{lem}[thm]{Lemma}

\theoremstyle{definition}
\newtheorem{defn}[thm]{Definition}

\newtheorem{exmp}[thm]{Example}
\newtheorem*{warn}{Warning}
\newtheorem{exmps}[thm]{Examples}

\theoremstyle{remark}
\newtheorem{rem}[thm]{Remark}
\newtheorem{rems}[thm]{Remarks}

\def\NN{\mathbf{N}}
\def\ZZ{\mathbf{Z}}

\def\RR{\mathbf{R}}

\def\defq{\stackrel{\text{\rm\tiny def}}{=}}
\newcommand{\overto}[1]{{\buildrel{#1}\over\longrightarrow}}
\newcommand{\overfrom}[1]{{\buildrel{#1}\over\longleftarrow}}

\newcommand{\setdef}[2]{ \left\{ {#1}\,\mid\, {#2} \right\} }
\newcommand{\ip}[2]{ \left\langle {#1} , {#2} \right\rangle }

\def\acts{\curvearrowright}

\def\normal{\triangleleft}
\newcommand{\ro}{\varrho}

\def\beq{\begin{equation}}
\def\eeq{\end{equation}}
\def\SL{{\bf\rm SL}}

\def\SO{{\bf\rm SO}}

\def\dist{{\rm dist}}
\def\diam{{\rm diam}}

\def\Ker{{\rm Ker}}
\def\Im{{\rm Im}}

\def\symdiff{{\,\triangle\,}}

\newcommand{\Rel}[1]{ \mathcal{#1} }

\def\Id{{\rm Id}}

\def\sH{{\mathcal{H}}}

\def\fintyp{{\mathscr{U}_{\rm fin}}}
\def\biinv{{\mathscr{G}_{\rm inv}}}
\def\cpt{{\mathscr{G}_{\rm cmp}}}
\def\dscr{{\mathscr{G}_{\rm dsc}}}

\def\Prob{{\sf Prob}}

\def\Aut{{\rm Aut}\,}
\def\Inn{{\rm Inn}\,}
\def\Out{{\rm Out}\,}
\def\Sym{{\rm Sym}}

\def\lcsc{{l.c.s.c. }} 

\title{On Popa's Cocycle Superrigidity Theorem}

\author[Alex Furman]{Alex Furman}
\address{University of Illinois at Chicago, Chicago, USA}
\email{furman@math.uic.edu}
\thanks{Supported in part by the NSF grants DMS-0094245, DMS-0604611, and BSF 2004345.}

\begin{document}

\begin{abstract}
These notes contain an Ergodic-theoretic account
of the Cocycle Superrigidity Theorem recently discovered by Sorin Popa.
We state and prove a relative version of the result,
discuss some applications to measurable equivalence relations,
and point out that Gaussian actions (of ``rigid'' groups) satisfy the 
assumptions of Popa's theorem. 
\end{abstract}

\maketitle

\section{Introduction and Statement of the Results}
A recent paper \cite{Popa:cocycle} of Sorin Popa contains a remarkable new Cocycle Superrigidity theorem. Here is a special case of this result
(see \cite[Theorem 3.1]{Popa:cocycle}, and Theorem~\ref{T:PopaCSR} below for
the general case):
\begin{thm}[{\bf Popa's Cocycle Superrigity, Special Case}]
\label{T:CSR-special}
Let $\Gamma$ be a discrete group with Kazhdan's property (T), 
let $\Gamma_{0}<\Gamma$ be an infinite index subgroup, $(X_{0},\mu_{0})$
an arbitrary probability space, and let 
$\Gamma\acts (X,\mu)=(X_{0},\mu_{0})^{\Gamma/\Gamma_{0}}$
be the corresponding generalized Bernoulli action.
Then for any discrete countable group $\Lambda$ and any measurable cocycle 
$\alpha:\Gamma\times X\to \Lambda$ there exist a homomorphism 
$\ro:\Gamma\to\Lambda$  and a measurable map $\phi:X\to\Lambda$
so that:   
$
	\ \alpha(g,x)=\phi(g.x)\ro(g)\phi(x)^{-1}.
$
\end{thm}

Note the unprecedented strength of this theorem:
there are no assumptions on the nature of the discrete target group $\Lambda$, 
no conditions are imposed on the cocycle $\alpha$, and the assumption on the
acting group $\Gamma$ is of a general and non-specific nature. 
At the same time the conclusion is unusually strong: 
the cocycle $\alpha$ is ``untwisted'' within the given discrete group $\Lambda$. 
The result however is specific to a particular class of actions -- Bernoulli actions,
or more generally \emph{malleable actions}.

This theorem continues a series of ground breaking results in the area of von Neumann Algebras
recently obtained by Sorin Popa  \cite{Popa:noncommB:06}, \cite{Popa:BettiAnnMath:06}, 
\cite{Popa:RigidityI}, \cite{Popa:RigidityII}, \cite{Popa:1cohom},  and by Popa and collaborators: 
Popa-Sasyk \cite{PopaSasyk}, Ioana-Peterson-Popa \cite{IoanaPetersonPopa:FreeProducts}, 
Popa-Vaes \cite{PopaVaes:06}.
The reader is also referred to Vaes' Seminare Bourbaki \cite{Vaes:afterPopa} for an overview
and more references.
The aims of these notes are:
\begin{enumerate}
\item 
	To give a rather short, self contained, purely Ergodic-theoretic proof of 
	Popa's Cocycle Superrigidity Theorem. We use this opportunity to give a
	relative version of the result (Theorem \ref{T:PopaCSR} below).
	Our proof follows Popa's general strategy, but implements some of the steps differently.
	The proof is essentially contained in sections \ref{S:Cohom} 
	and \ref{S:rigidity-deformation} below.
\item
	To point out that the class of actions for which Popa proved his 
	Cocycle Superrigidity property (weakly mixing \emph{malleable} 
	actions) extends from the class of (generalized) Bernoulli actions to 
	a larger class of all weakly mixing Gaussian actions.
\item
	To discuss the relationship between the Cocycle Superrigidity results
	of Zimmer and Popa, showing in particular that algebraic actions 
	are not quotients of malleable actions.
\item
	To discuss some applications of Popa's Cocycle Superrigidity Theorem to Ergodic theory,
	in particular to ${\rm II}_{1}$ equivalence relations.
\end{enumerate}
In these notes we emphasize the Ergodic theoretic point of view, trying to compliment the original 
Operator Algebra framework of Popa's \cite{Popa:cocycle} (also Vaes' \cite{Vaes:afterPopa}). 
As a result we do not discuss at all the results on von Neumann equivalence
of group actions, and refer the reader to the above mentioned papers 
(see also \cite{PopaVaes:06}).

\subsection{The Statement of Popa's Cocycle Superrigidity Theorem}
\label{SS:mainthm}

Theorem \ref{T:CSR-special} is only a special case of Popa's Cocycle Superrigidity theorem.
Before formulating the more general statement \ref{T:PopaCSR} below it is necessary to 
recall/define some notions which we shall use hereafter:
\begin{description}
\item[w-normal and wq-normal subgroups]
The following notions are generalizations of normality, and subnormality.
A closed subgroup $H<G$ is called \emph{weakly normal} (\emph{w-normal}),
if there exists an increasing family $H=H_{0}\normal H_{1}\normal \cdots H_{\eta}=G$ 
of closed subgroups $\{H_{i}\}_{0\le i\le \eta}$, well ordered by inclusion, and such 
that $\left(\bigcup_{i<j} H_{i}\right)\normal H_{j}$ for each ordinal $0\le j\le \eta$.
A further weakening of this notion is: 
$H<G$ is \emph{weakly quasi-normal} (\emph{wq-normal}), 
if there exists an increasing family $H=H_{0}< H_{1}\normal \cdots H_{\eta}=G$ 
of closed subgroups $\{H_{i}\}_{0\le i\le \eta}$ well ordered by inclusion, such 
that denoting $H'_{j}\defq\left(\bigcup_{i<j} H_{i}\right)$ for each ordinal $0\le j\le \eta$ 
the group $H_{j}$ is generated by the set 
$\setdef{g\in G}{g^{-1}H'_{j}g\cap H'_{j}\ \text{is not precompact in } G}$. 

\item[groups of finite type]
A topological group is said to be of \emph{finite type}, or to belong to class $\fintyp$,
if it can be embedded as a closed subgroup of the unitary group of a finite 
von-Neumann algebra (i.e. a von Neumann algebra with a faithful normal final trace).
Class $\fintyp$ contains the class $\dscr$ of all discrete countable groups, 
and the class $\cpt$ of all second countable compact groups.
$\fintyp$ also contains such groups as the inner automorphism group 
$\Inn(\Rel{R})$ (a.k.a. the full group)
of a II$_{1}$ countable relation $\Rel{R}$. 
All groups in $\fintyp$ admit a complete metric which is bi-invariant;
hence semi-simple Lie groups are not in $\fintyp$ 
(see \ref{SS:Targetgroups} for further discussion).  

\item[(relative) property (T)]
A locally compact second countable (\lcsc) group $G$ has Kazhdan's property (T) 
if every unitary $G$-representation 
which almost has invariant vectors, has non-trivial invariant vectors.
A closed subgroup $H<G$ in a \lcsc group $G$ is said to have \emph{relative property (T)} in $G$ 
if every unitary $G$-representation which almost has $G$-invariant vectors, 
has non-trivial $H$-invariant vectors.

\item[$L$-Cocycle-Superrigid actions (extensions)]
Let $G\acts (X,\mu)$ be an ergodic probability measure preserving
(p.m.p.) action of some \lcsc group $G$,
and $L$ be some Polish group. 
We shall say that the action $G\acts (X,\mu)$ is $L$-\emph{Cocycle-Superrigid}
if for every measurable cocycle $\alpha:G\times X\to L$ there 
exists a measurable map $\phi:X\to L$ and a homomorphism $\ro:G\to L$
so that for each $g\in G$: 
\[
	\alpha(g,x)=\phi(g.x) \ro(g)\,\phi(x)^{-1}
\]
for $\mu$-a.e. $x\in X$.
More generally, if $p:(X,\mu)\to (Y,\nu)$ is a measurable equivariant quotient map of p.m.p.
$G$-actions, we shall say that $G\acts (X,\mu)$ is $L$-\emph{Cocycle-Superrigid relatively to} 
$(Y,\nu)$ if every measurable cocycle $\alpha:G\times X\to L$ can be written as
\[
	\alpha(g,x)=\phi(g.x) \ro(g,p(x))\,\phi(x)^{-1}
\]
where $\phi:X\to L$ is a measurable map and $\ro:G\times Y\to L$
is a measurable cocycle.
If an action $G\acts (X,\mu)$ is  $L$-Cocycle-Superrigid for all groups $L$ in some class
$\mathscr{C}$ (e.g. $\fintyp$)
we shall say that  $G\acts (X,\mu)$ is $\mathscr{C}$-\emph{Cocycle-Superrigid}.
Similarly one defines the relative notion.

\item[weakly mixing actions (extensions)]
A p.m.p action $G\acts (X,\mu)$ is said to be \emph{weakly mixing}
if the diagonal action of $G$ on $(X\times X,\mu\times\mu)$ is ergodic.
If $G\acts (X,\mu)$ has a quotient p.m.p. action $G\acts (Y,\nu)$ 
one can form a \emph{fibered product space} 
$X\times_{Y} X\defq\setdef{(x_1,x_2)\in X\times X}{p(x_1)=p(x_2)}$
on which $G$ acts diagonally preserving a certain probability measure $\mu\times_{\nu}\mu$,
associated to the disintegration of $p:\mu\mapsto \nu$. 
If the latter action is ergodic, the original action $G\acts (X,\mu)$ is 
said to be \emph{weakly mixing relative to} the quotient $G\acts (Y,\nu)$;
weak mixing is equivalent to weak mixing relative to the trivial action on a point
(see \ref{SS:relwm} for more details). 

\item[malleable actions (extensions)]
Given a probability space $(Z,\zeta)$ denote by $\Aut(Z,\zeta)$ the group
of all measure space automorphisms of $(Z,\zeta)$ considered modulo null sets
and endowed with the \emph{weak topology}
which makes it into a Polish group (see \ref{SS:topologies}).
A p.m.p action $G\acts (X,\mu)$  is \emph{malleable} if the connected component of the 
identity of the centralizer $\Aut(X\times X,\mu\times\mu)^{G}$ of the diagonal $G$-action
on $(X\times X,\mu\times\mu)$  contains the flip $F:(x,y)\mapsto(y,x)$ or, 
an element of the form 
\[
	\qquad F\circ(T\times S):(x,y)\mapsto (S(y),T(x))\quad\text{where}
	\quad T,S\in\Aut(X,\mu)^{G}.
\]
More generally, if $(X,\mu)\overto{p} (Y,\nu)$ is a morphism of p.m.p. $G$-actions
we shall say that $G\acts (X,\mu)$ is \emph{malleable relative to} $(Y,\nu)$
if the flip $F:(x_{1},x_{2})\mapsto (x_{2},x_{1})$,
or a transformation of the form $F\circ (T\times S)$ with $T,S\in\Aut(X,\mu)^{G}$,
$p\circ T=p=p\circ S$, lies in the connected component of the identity 
in the centralizer  of $G$ in $\Aut(X\times_{Y} X,\mu\times_{\nu}\mu)$.
\footnote{We use a variation of the original definitions (cf. \cite[4.3]{Popa:cocycle}),
using connectivity rather than path connectivity and allowing for $S,T$.}

\item[generalized Bernoulli actions]	
Given a probability space $(X_{0},\mu_{0})$, 	
an infinite countable set $I$, and a permutation action $\sigma:G\to \Sym(I)$,
there is a p.m.p. $G$-action on the product space 
$(X,\mu)=(X_{0},\mu_{0})^{I}=\prod_{i\in I} (X_{0},\mu_{0})_{i}$ 
by  permutation of the coordinates: $(g.x)_{i}=(x)_{\sigma(g^{-1})(i)}$. 
This action $G\acts (X,\mu)$ is ergodic iff the $G$-action on the index set $I$ 
has no finite orbits.
The classical (two-sided) Bernoulli shift corresponds to $G=\ZZ$ acting on itself.

\item[Gaussian actions]
are p.m.p. actions $G\acts (X,\mu)$ constructed out of unitary, or rather orthogonal, representations
$\pi:G\to O(\sH)$ with the $G$-representation on $L^{2}(X,\mu)$ being isomorphic to the direct sum
$\oplus_{n=0}^{\infty} S^{n}\pi$ of the symmetric powers of $\pi$, which in particular contains $\pi$
itself as a subrepresentation. 
For Gaussian actions ergodicity is equivalent to weak mixing and occurs iff $\pi$ is weakly mixing 
(i.e., does not contain a finite dimensional subrepresentation).
We note that if $G$ acts by permutations on a countable set $I$, then the Bernoulli
$G$-action on $([0,1],m)^{I}$ is isomorphic to the Gaussian
action corresponding to the $G$-representation on $\ell^{2}(I)$;
hence the class of Gaussian actions contains that of generalized Bernoulli 
actions with the non atomic base space.
\end{description}
The notions of \emph{w-normality} and \emph{wq-normality}, 
the class $\fintyp$ of groups, and the concept
 of \emph{malleability} (and other variants of these notions) were introduced 
 and studied by Sorin Popa. 
The notion of malleability is of particular importance in his work 
and in our discussion of it. 
Generalized Bernoulli actions (with a non-atomic base space) 
were the main example of malleable actions in \cite{Popa:cocycle}. 
Here we observe that the class of malleable actions contains all  Gaussian actions.
\begin{thm}
\label{T:malleable}
All Gaussian actions, including generalized Bernoulli actions 
with a non-atomic base space, of any \lcsc group are malleable. 
\end{thm}
Note that the notion of malleability itself does not involve any ergodicity assumptions,
and so in the above statement no assumptions are imposed on the Bernoulli or Gaussian 
actions. Yet in what follows we shall need to require (weak) mixing of the action of the group $G$,
or even of a certain subgroup $H<G$.
For a Bernoulli action $G\acts (X_{0},\mu_{0})^{I}$ the subgroup $H$ would be weakly mixing
iff its action on the index set $I$ has no finite orbits (we always assume that the base space
$(X_{0},\mu_{0})$ is not a point, but allow it to have atoms).
For Gaussian actions arising from an orthogonal $G$-representation $\pi$, 
the action of a subgroup $H<G$ is weakly mixing iff the restriction $\pi|_{H}$ has 
no finite dimensional subrepresentations.
We are ready to state the general form of Popa's Cocycle Superrigidity Theorem 
(cf. \cite[Theorem 3.1]{Popa:cocycle}). 
\begin{thm} [{\bf Popa's Cocycle Superrigidity (with the relative version)}]
\label{T:PopaCSR}
Let $G$ be a \lcsc group with a closed subgroup $H<G$ which has relative property (T) in $G$,
and let $G\acts (X,\mu)$ be an ergodic p.m.p. action (with a quotient $(Y,\nu)$).
Suppose that:
\begin{itemize}
\item[(a)]
    $G\acts (X,\mu)$ is malleable (resp. malleable relatively to $(Y,\nu)$), and either
\item[(b)]
    $H$ is w-normal in $G$ and $H\acts (X,\mu)$ is weakly mixing 
    (resp. weakly mixing relatively to $(Y,\nu)$), or
\item[(b')]
    $H$ is wq-normal in $G$ and $G\acts (X,\mu)$ is mixing 
    (resp. mixing relatively to $(Y,\nu)$\footnote{The actual condition is
    that for any non pre-compact subgroup $G'<G$ the action
    $G'\acts (X,\mu)$ is weakly mixing (resp. relatively to $(Y,\nu)$).}).
\end{itemize} 
Then the $G$-action on $(X,\mu)$ is $\fintyp$-Cocycle-Superrigid 
(resp.  relatively to $(Y,\nu)$).
\end{thm}
For a non-atomic probability space $(X_{0},\mu_{0})$ the combination 
of the above Theorem~\ref{T:PopaCSR} and the malleability 
(\ref{T:malleable}) imply $\fintyp$-Cocycle Superrigidity of Bernoulli actions 
$\Gamma\acts (X_{0},\mu_{0})^{\Gamma/\Gamma_{0}}$
as in Theorem~\ref{T:CSR-special}.
To deduce $\fintyp$-Cocycle Superrigidity for Bernoulli actions 
where the base space $(X_{0},\mu_{0})$ has atoms, Popa uses 
the following argument (cf. \cite[Lemma 2.11]{Popa:cocycle},
see Lemma~\ref{L:basic} below):
\begin{prop}[{\bf Weakly Mixing Extensions}]
\label{P:CSRrelwm}
Let $G$ be a \lcsc group, $(X',\mu')\to (X,\mu)\to (Y,\nu)$ be p.m.p extensions
with $(X',\mu')\to (X,\mu)$ being relatively weakly mixing. 
If  $G\acts (X',\mu')$ is $\fintyp$-Cocycle-Superrigid (relatively to $(Y,\nu)$), 
then the same applies to the intermediate action $G\acts (X,\mu)$.
\end{prop}
Indeed, any probability space $(X_{0},\mu_{0})$ is a quotient of \emph{the} 
non-atomic probability space $([0,1],m)$ where $m$ is the Lebesgue measure; 
hence $(X,\mu)=(X_{0},\mu_{0})^{I}$ is a $G$-quotient of $(X',\mu')=([0,1],m)^{I}$,
which is easily seen to be relatively weakly mixing.

\medskip

Skew-products provide examples of \emph{relative} $\fintyp$-Cocycle Superrigid actions.
\begin{exmp}[{\bf Skew products}]
Let $H<G$ be a w-normal subgroup with relative property (T) with a p.m.p. action
$G\acts (Y,\nu)$ and a measurable cocycle $\sigma:G\times Y\to V$ with values
in a \lcsc group $V$, and let $V\acts (Z,\zeta)$ be a malleable p.m.p. action,
e.g. generalized Bernoulli or Gaussian action.
Let $G$ act on $(X,\mu)=(Y\times Z,\nu\times\zeta)$ by $g:(x,z)\mapsto (g.x,\sigma(g,x).z)$.
If the $H$-action on $(X,\mu)$ is weakly mixing relative to $(Y,\nu)$, i.e.,
if $H$ acts ergodically on $Y\times Z\times Z$, then $G\acts (X,\mu)$ is $\fintyp$-Cocycle Superrigid
relatively to $(Y,\nu)$.
\end{exmp}
First we discuss some applications of the ``absolute case'' of the theorem to
equivalence relations, presenting the Ergodic theoretic rather than Operator Algebra point of view.
For very interesting further applications see the recent paper of  
Popa - Vaes \cite{PopaVaes:06}.

\bigskip

\subsection{Applications to II$_1$ Equivalence Relations}
\label{SS:OE}
Let $\Gamma\acts (X,\mu)$ be an ergodic action of a countable group $\Gamma$.
The ``orbit structure'' of this action is captured by the equivalence relation
\[
	\Rel{R}_{X,\Gamma}=\setdef{(x,x')\in X\times X}{\Gamma. x=\Gamma. x'}.
\]
Two actions $\Gamma\acts (X,\mu)$ and $\Lambda\acts (Y,\nu)$ are said to be
\emph{Orbit Equivalent} (OE) if their orbit relations are isomorphic, 
where isomorphism of equivalence relations, is a measure space isomorphism
$\theta:X'\to Y'$ between conull subsets of $X$ and $Y$ identifying the restrictions 
of the corresponding equivalence relations to these subsets.
Isomorphism of equivalence relations and OE of actions 
can be weakened by allowing $X'$ and $Y'$ to be arbitrary
measurable subsets of \emph{positive} (rather than full) measure.
This is  \emph{weakly isomorphism} of relations, 
and \emph{weakly Orbit Equivalent} (or \emph{stably Orbit Equivalent}) of actions;
if the original actions are of type ${\rm II}_{1}$ the ratio $\nu(Y')/\mu(X')$ is called the
\emph{compression constant} of $\theta$ (\cite{Furman:OE:99}). 
Zimmer's Cocycle Superrigidity Theorem \cite{Zimmer:cocyclesuper:80}
paved the way to many very strong OE rigidity results in ergodic theory for higher rank lattices
and other similar groups (see Zimmer's monograph \cite{Zimmer:book:84} and references therein, 
Furman \cite{Furman:ME:99}, \cite{Furman:OE:99}, \cite{Furman:Outer:05},
Monod - Shalom \cite{MonodShalom:OE:05}, Hjorth - Kechris \cite{HjorthKechris:MAMS:05}),
also (Gaboriau \cite{Gaboriau:Inventcost:00}, \cite{Gaboriau:L2IHES:02}, 
and Shalom \cite{Shalom:MGT:05} for a survey).

Popa's Cocycle Superrigidity Theorem \ref{T:PopaCSR} allows to study not only (weak) OE 
isomorphisms, but also (weak) OE ``imbeddings'' and other ``morphisms'' between equivalence relations. 
By a \emph{morphism} $\theta:(X,\mu,\Rel{R})\to (Y,\nu,\Rel{S})$ between ergodic countable equivalence relations 
we mean a measurable map $\theta:X'\to Y'$ between conull spaces $X'\subset X$, $Y'\subset Y$ with 
$\theta_{*}\mu|_{X'}\sim\nu|_{Y'}$ and 
$\theta\times\theta(\Rel{R}|_{X'\times X'})\subset\Rel{S}|_{Y'\times Y'}$.
The \emph{kernel} of the morphism $\theta$ is the relation on $(X,\mu)$ given by
\[
	\Ker(\theta)=\setdef{(x_{1},x_{2})\in\Rel{R}}{\theta(x_{1})=\theta(x_{2})}
\] 
(note that $\Ker(\theta)$ it is not ergodic if $Y$ is a non-trivial space);
the \emph{image} $\Im(\theta)$ of $\theta$ is the obvious subrelation 
of $\Rel{S}$ on $(Y,\nu)$.
\begin{defn}[{\bf Relation Morphisms}]
Let $\theta$ be a weak relation morphism from an ergodic 
equivalence relation $\Rel{R}$ on $(X,\mu)$ to $\Rel{S}$ on $(Y,\nu)$. 
We shall say that $\theta$ is an \emph{injective relation morphism} 
if $\Ker(\theta)$ is trivial, 
that $\theta$ is a \emph{surjective relation morphism} if $\Im(\theta)=\Rel{S}$, 
and that $\theta$ is a \emph{bijective relation morphism} 
if $\Ker(\theta)$ is trivial and $\Im(\theta)=\Rel{S}$.
The notion of \emph{weakly injective/surjective/bijective relation morphism}
corresponds to $X'$ and $Y'$ being positive measure subsets, rather than conull ones.
In the ${\rm II}_{1}$ context weak morphisms come equipped with a compression 
constant\footnote{In \cite{Popa:cocycle} a bijective relation morphism is called 
\emph{local Orbit Equivalence}.}.
\end{defn}

\begin{rem}
Note that a bijective relation morphism is not necessarily an isomorphism, 
because $\theta:X\to Y$ need not be invertible.
For example, any $\Gamma$-equivariant quotient map $\theta:(X,\mu)\to(Y,\nu)$ 
of free ergodic \emph{actions} of some countable group 
$\Gamma$ defines a bijective relation morphism of the corresponding orbit relations 
$(X,\mu,\Rel{R}_{X,\Gamma})\to(Y,\nu,\Rel{R}_{Y,\Gamma})$.
\end{rem}
The following theorem summarizes some of the  consequences of Theorem~\ref{T:PopaCSR}
to orbit relations;
it is parallel to and somewhat more detailed than 
Popa's \cite[Theorems 0.3, 0.4, 5.6--5.8]{Popa:cocycle}.
We denote by  $\dscr\subset\fintyp$ the class of all discrete countable groups.
\begin{thm}[{\bf Superrigidity for Orbit Relations}]
\label{T:superOE}
Let $\Gamma\acts (X,\mu)$ be a $\dscr$-Cocycle Superrigid p.m.p. ergodic action of a 
countable group $\Gamma$, and $\Lambda\acts (Y,\nu)$ be an ergodic 
measure preserving  essentially free action of some countable group $\Lambda$
on a (possibly infinite) measure space $(Y,\nu)$ and let 
\[
	X\supset X'\overto{\theta} Y'\subset Y
\] 
be a weak morphism of the orbit relations 
$\Rel{R}_{X,\Gamma}$ and $\Rel{R}_{Y,\Lambda}$. Then there exist: 
\begin{enumerate}
\item 
	a homomorphism $\ro:\Gamma\to\Lambda$ (we denote $\Gamma_{0}=\Ker(\ro)$ and 
	$\Gamma_{1}=\ro(\Gamma)<\Lambda$);
\item 
	a measurable $\Gamma_{1}$-invariant subset $Y_{1}\subset Y$
	with $0<\nu(Y)<\infty$ (we denote $\nu_{1}= \frac{1}{\nu(Y_{1})}\cdot\nu_{|Y_{1}}$)
	with the p.m.p. action $\Gamma_{1}\acts (Y_{1},\nu_{1})$ being ergodic;
\item 
	a measurable map $T:X\to Y_{1}$ with $T_{*}\mu=\nu_{1}$ and 
	\[
		T(g.x)=\ro(g).T(x)\qquad (g\in\Gamma)
	\] 
\end{enumerate}
The map $T$ is a composition
	$
		T:(X,\mu)\overto{erg}(X_{1},\mu_{1})\overto{T_{1}}(Y_{1},\nu_{1}),
	$
where $(X_{1},\mu_{1})$ is the space of $\Gamma_{0}$-ergodic components equipped
with the natural action of $\Gamma_{1}\cong \Gamma/\Gamma_{0}$, and 
$T_{1}:(X_{1},\mu_{1})\to (Y_{1},\nu_{1})$ is a quotient map of $\Gamma_{1}$-actions.

Furthermore, with the notations above:
\begin{itemize}
\item[{\rm (i)}]
	If $\Gamma\acts (X,\mu)$ is essentially free, then 
	$\theta$ is a weakly injective morphism iff $\Gamma_{0}=\Ker(\ro)$ is finite;
	in this case $(X,\mu)$ is a finite extension of $(X_{1},\mu_{1})$.
\item[{\rm (ii)}]
	If $\nu(Y)<\infty$ and $\theta$ is a weakly surjective morphism, then 
	$\Gamma_{1}=\ro(\Gamma)$ is of finite index in $\Lambda$
	and, assuming $\Gamma\acts(X,\mu)$ is aperiodic, 
	$\nu(Y)/\nu(Y_{1})$ is an integer dividing $[\Lambda:\Gamma_{1}]$.
\item[{\rm (iii)}]
	If $\nu(Y)=1$, $\Gamma\acts (X,\mu)$ is essentially free and aperiodic,
	and $\theta$ is a weakly bijective morphism (in particular, if $\theta$
	is a weak isomorphism) of $\Rel{R}_{X,\Gamma}$ to $\Rel{R}_{Y,\Lambda}$, 
	then both $|\Gamma_{0}|$ and
	$[\Lambda:\Gamma_{1}]$ are finite, the compression constant 
	is rational and is given by 
	\[
		c(\theta)=\frac{\nu(Y_{1})}{\nu(Y)}\cdot\frac{[\Lambda:\Gamma_{1}]}{|\Gamma_{0}|}.
	\] 
	Such a weakly bijective morphism (or weak isomorphism) $\theta$ extends to a bijective 
	morphism (resp. an isomorphism) $\tilde\theta:X\to Y$ iff $c(\theta)=1$.
\end{itemize}
\end{thm}
\begin{rems}
\begin{enumerate}
\item
An action $\Gamma\acts (X,\mu)$ is called aperiodic, if it has no non-trivial finite quotients;
equivalently, if the restriction of the action $\Gamma'\acts (X,\mu)$ to any finite index
subgroup $\Gamma'<\Gamma$ is ergodic. Weakly mixing actions are aperiodic.
\item
Note that often $\Gamma_{0}$ is forced to be finite, even in the context of general weak morphisms 
of orbit relations (coming from $\dscr$-Cocycle Superrigid action to a free action of some 
countable group).
Finiteness of $\Gamma_{0}$ is guaranteed if every infinite normal subgroup of 
$\Gamma$ acts ergodically on $(X,\mu)$; which is indeed the case if $\Gamma\acts(X,\mu)$ is mixing,
or if $\Gamma$ has no proper infinite normal subgroups (e.g. if $\Gamma$ is an irreducible
lattice in a higher rank Lie group, see \cite[Ch VIII]{Margulis:book:91}).
\item
In (ii) one cannot reverse the implication. If $\Gamma<\Lambda$ is a finite index
inclusion of countable groups with property (T) so that $\Gamma$ has no epimorphism
onto $\Lambda$, and $(X,\mu)=(X_{0},\mu_{0})^{\Lambda}$ with the Bernoulli action of 
the groups $\Gamma<\Lambda$,
then the identity map $\theta(x)=x$ is a relation morphism from $\Rel{R}_{X,\Gamma}$ to
$\Rel{R}_{X,\Lambda}$ which is not weakly surjective.
In fact there does not exist weakly surjective morphisms from 
$\Rel{R}_{X,\Gamma}$ to $\Rel{R}_{X,\Lambda}$.
\end{enumerate}
\end{rems}
The rigidity with respect to weak self isomorphisms as in (iii) 
provides the basic tool for the computation of certain invariants such as the fundamental group, 
first cohomology, the outer automorphism group of the orbit relations of 
$\fintyp$-Cocycle Superrigid actions.
It can also be used to produce equivalence relations which cannot be 
generated by an essentially free action 
of any group (see \cite{PopaSasyk}, \cite{PopaVaes:06},
for actions of higher rank lattices see \cite{Furman:OE:99}, \cite{Furman:Outer:05}). 
It seems, however, that the full strength of Popa's Cocycle Superrigidity Theorem, manifested in the
first part of Theorem~\ref{T:superOE}, should yield qualitatively new rigidity phenomena related 
to morphisms of relations.

\subsection{Comparison with Zimmer's Cocycle Superrigidity}
\label{SS:Zimmer}
It is natural to compare Popa's Theorem with Zimmer's Cocycle Superrigidity. 
Zimmer's result (see \cite[Theorem 5.2.5]{Zimmer:book:84}),
generalizing Margulis' Superrigidity (see \cite[Ch VII]{Margulis:book:91}), is a theorem
about untwisting of measurable cocycles $\alpha:G\times X\to H$ 
over ergodic (irreducible) actions of (semi-)simple algebraic groups $G$
and taking values in semi-simple algebraic groups $H$
(the result is subject to the assumption that $\alpha$ is Zariski dense and is not ``precompact'').
Zimmer's cocycle superrigidity has a wide variety of applications with two most prominent 
areas being:

(a) smooth volume preserving actions of higher rank groups/lattices
on manifolds (here the cocycle is the derivative cocycle),

(b) orbit relations of groups actions in Ergodic theory (using the ``rearrangement'' cocycle). 

The comparison with Popa's result lies in the latter area.
In this context Zimmer's Cocycle Superrigidity Theorem is usually applied
to cocycles $\alpha:\Gamma\times X\to \Lambda$, where $\Gamma<G$ is a lattice
in a simple Lie group, and $\Lambda$ is a lattice or a more general subgroup
in a (semi)simple Lie group $H$.
Assuming the cocycle satisfies the assumptions of the theorem (which is the case 
if $\alpha$ comes from an OE, or weak OE of free p.m.p. actions of 
$\Gamma\acts (X,\mu)$ and $\Lambda\acts (X',\mu')$) the conclusion is that: 
	there exists a homomorphism (local isomorphism) $\ro:G\to H$
	and a measurable map $\phi:X\to H$ so that
	$\alpha(\gamma,x)=\phi(\gamma.x)\,\ro(\gamma)\,\phi(x)^{-1}$.
In comparison with Popa's result note the following points:
\begin{itemize}
\item[(i)]
	there is no claim that $\ro(\Gamma)$ is a subgroup of $\Lambda<H$;
\item[(ii)]
	even if this is the case, say if $G=H$, $\Gamma=\Lambda$ and $\ro$ is the identity,
	one still cannot expect the ``untwisting'' map $\phi(x)$ to take values in $\Lambda<H$;
\item[(iii)]
	at the same time there is no specific assumptions on the action $\Gamma\acts (X,\mu)$
	beyond ergodicity (or irreducibility if $G$ is semi-simple).
\end{itemize}	
The point (ii) is extensively studied in \cite{Furman:Outer:05}. The simplest example
of a cocycle $\Gamma\times X\to \Gamma<G$ which is cohomologous to the identity 
homomorphism as a cocycle into $G$ but not as a cocycle into $\Gamma$, appears in
the left translation $\Gamma$-action on $X=G/\Gamma$.
The cocycle is $\alpha=\ro|_{\Gamma\times G/\Gamma}$ --
the restriction of the \emph{canonical cocycle} $\ro:G\times G/\Gamma\to\Gamma$
defined by $\ro(g,x)=f(g.x)\,g\,f(x)^{-1}$,
where $f:G/\Gamma\to G$ is a Borel cross-section of the projection $G\to G/\Gamma$.
So Zimmer's Cocycle Superrigidity does not, and in general cannot, untwist cocycles within
a discrete group, even if the target group is $\Gamma$ itself.
This is, of course, in a sharp contrast to Popa's result.

Another feature of Popa's Cocycle Superrigidity is that it has no assumptions on
the target group besides it being discrete, or more generally, in class $\fintyp$.
The main point of \cite{Furman:ME:99} and \cite{Furman:OE:99} was to prove
a result of this type for cocycles $\alpha:\Gamma\times X\to\Lambda$ 
over a general ergodic action $\Gamma\acts (X,\mu)$ of a higher rank lattice $\Gamma$,
with $\Lambda$ being an arbitrary discrete countable group.
These results, however, require the following condition on the cocycle $\alpha$: 
\begin{defn}
\label{D:finvol}
Let $\Gamma$ and $\Lambda$ be discrete countable groups, $\Gamma\acts (X,\mu)$
be a p.m.p. action and $\alpha:\Gamma\times X\to\Lambda$ be a measurable cocycle.
We shall say that $\alpha$ is a \emph{Measure Equivalence} cocycle (ME-cocycle), 
or \emph{finite covolume} cocycle, if the $\Gamma$-action on the infinite Lebesgue space
\[
	(\tilde{X},\tilde{\mu})=(X\times\Lambda,\mu\times m_{\Lambda})
	\qquad
	g:(x,\ell)\mapsto (g.x, \alpha(g,x)\ell)
	\qquad(g\in\Gamma)
\] 
is essentially free, and has a Borel fundamental domain of finite $\tilde{\mu}$-measure. 
\end{defn}
\begin{exmp}
Cocycles corresponding to (weak) Orbit Equivalence of essentially
free actions have finite covolume. The same can be shown for cocycles
corresponding to (weakly) bijective morphisms of orbit relations of essentially free p.m.p. actions.
As the name suggests cocycles associated to any Measure Equivalence coupling
between two countable groups are ME-cocycles (see \cite{Furman:OE:99}).
\end{exmp} 
\begin{thm}[{\bf ME-Cocycle Superrigidity}]
\label{T:ME+}
Let $k$ be a local field, $\mathbf{G}$ be a $k$-connected simple algebraic
group defined over $k$ with ${\rm rank}_{k}(\mathbf{G})\ge 2$, $G=\mathbf{G}(k)$ 
the \lcsc group of the $k$-points of $\mathbf{G}$.
Denote $\bar{G}=\Aut(G/Z(G))$ and let $g\mapsto\bar{g}$ be the homomorphism
$G\to\bar{G}$ with finite kernel $Z(G)$ and finite cokernel $\Out(G/Z(G))$. 

Let $\Gamma<G$ be a lattice and $\Gamma\acts (X,\mu)$ be an ergodic 
aperiodic p.m.p. action, and let $ A=\{\alpha_{i}\}$ be a family of measurable 
cocycles $\alpha_{i}:\Gamma\times (X,\mu)\to\Lambda_{i}$ of finite covolume 
taking values in arbitrary countable groups $\Lambda_{i}$.

Then there exists a measurable $\Gamma$-equivariant quotient 
$(X,\mu)\overto{p} (Y_{ A},\nu_{ A})$ with the following properties:
\begin{enumerate}
\item
	for each cocycle $\alpha_{i}:\Gamma\times (X,\mu)\to\Lambda_{i}$ in $ A$
	there exist: a short exact sequence $N_{i}\to \Lambda_{i}\to\bar{\Lambda}_{i}$
	with $N_{i}$ being finite, and $\bar{\Lambda}_{i}$ being a lattice in $\bar{G}$;
	a measurable map $\phi_{i}:X\to\bar{\Lambda}_{i}$, and a measurable cocycle
	$\ro_{i}:\Gamma\times Y_{ A}\to\bar{\Lambda}_{i}$
	so that, denoting $\bar\alpha_{i}:\Gamma\times X\to\Lambda_{i}\to\bar\Lambda_{i}$,
	we have 
	\[
		\bar\alpha_{i}(\gamma,x)=\phi_{i}(\gamma.x)\,\ro_{i}(\gamma,p(x))\,\phi_{i}(x)^{-1}.
	\] 
\item
	$(Y_{ A},\nu_{ A})$ is the minimal quotient with the above properties,
	i.e., if all $\alpha\in A$ reduce to a quotient $(Y',\nu')$ then
	$(Y',\nu')$ has $(Y_{ A},\nu_{ A})$ as a quotient.
\end{enumerate}
In the case of a one point set $ A=\{\alpha\}$, where
$\alpha:\Gamma\times X\to\Lambda$ is a finite covolume cocycle,
the quotient $(Y_{\{\alpha\}},\nu_{\{\alpha\}})$ is either
\begin{itemize}
\item[{\rm (i)}]
 	trivial, in which case $\ro_{\alpha}:\Gamma\to\bar\Lambda$ is an inner homomorphism 
	$\ro_{\alpha}(\gamma)=\bar{g}_{0}\bar\gamma\bar{g}_{0}^{-1}$ 
	by some $\bar{g}_{0}\in\bar{G}$ (in particular, 
	$\bar{g}_{0}\bar\Gamma\bar{g}_{0}^{-1}<\bar\Lambda$), or 
\item[{\rm (ii)}]
	is $(\bar{G}/\bar\Lambda,m_{\bar{G}/\bar\Lambda})$
	and  $\ro_{\alpha}:\Gamma\times\bar{G}/\bar\Lambda\to\bar\Lambda$ 
	is the restriction of the standard cocycle. 
\end{itemize} 
In the general case $(Y_{ A},\nu_{ A})$ is the join  
$\bigvee_{\alpha\in A}(Y_{\{\alpha\}},\nu_{\{\alpha\}})$.
For a finite set $A=\{\alpha_{1},\dots,\alpha_{n}\}$ of cocycles 
the quotient $\Gamma\acts (Y_{A},\nu_{A})$ is isomorphic to
the diagonal $\Gamma$-action on $(\prod_{j=1}^{k}\bar{G})/\Delta$,
where $\Delta$ is a  lattice in $\prod_{j=1}^{k}\bar{G}$ containing a product
$\prod_{j=1}^{k} \Lambda_{i_{j}}$ of lattices in the factors.
\end{thm}
\begin{rems}
\begin{enumerate}
\item
Assuming that the target group $\Lambda$ has only infinite conjugacy classes 
(ICC) one has $N_{i}=\{e\}$, $\bar{\Lambda}_{i}=\Lambda_{i}$, 
$\bar\alpha_{i}=\alpha_{i}$.
\item
We recall that  many aperiodic actions $\Gamma\acts(X,\mu)$ do not have any measurable 
quotients of the form $\bar{G}/\Lambda$ (see \cite{Furman:OE:99}); 
for such actions case (a) is ruled out, and so in this case we get an absolute 
ME-Cocycle Superrigidity result (for general countable ICC targets). 
\item
The above Theorem is an extension of \cite[Theorems B,C]{Furman:OE:99},
the latter corresponds to the single ME-cocycle case $(Y_{\{\alpha\}},\nu_{\{\alpha\}})$. 
We have also used this opportunity to extend the framework from higher rank \emph{real} 
simple Lie groups to higher rank simple algebraic groups over general fields
(this extension does not require any additional work, the original proofs in
\cite{Furman:ME:99}, \cite{Furman:OE:99} applied verbatim to simple algebraic groups
over general local field).
\item
The definition~\ref{D:finvol} applies not only to discrete groups but to all
\emph{unimodular} \lcsc groups $\Lambda$.
Theorem~\ref{T:ME+} applies to this more general setting with the following adjustments: 
in (1) the kernel $N_{i}$ is compact, rather than just finite, and $\bar\Lambda_{i}$ 
in addition to being a lattice in $\bar{G}$ may also be $G/Z(G)$, $\bar{G}$,
or any intermediate closed group;
statements (i) and (ii) are claimed only if $\bar{\Lambda}$ is discrete (i.e., is a lattice
in $\bar{G}$).
\item
The results in \cite{Furman:cohom} suggest that in many situations it is possible to define
a canonical \emph{Cocycle Superrigid quotient} $(Y_{A},\nu_{A})$ for a given p.m.p. action 
$G\acts (X,\mu)$  and a set $A$ of measurable
cocycles $\alpha:G\times X\to L_{\alpha}$ with target groups $L_{\alpha}\in\biinv$.
\end{enumerate}
\end{rems}

\bigskip

\subsection{Standard algebraic actions vs. Malleable actions }
\label{SS:slg}
Continuing the line of comparison between Popa's and Zimmer's Cocycle Superrigidty
Theorems it is natural to ask whether $\fintyp$-Cocycle 
Superrigid actions $\Gamma\acts (X,\mu)$ of a higher rank lattice $\Gamma<G$ 
admit an algebraic action $\Gamma\acts G/\Lambda$ as a measurable quotient. 
We answer this question (posed by Popa) in the negative. 
In fact, a much more general class of algebraic actions cannot appear as a quotient of 
a $\fintyp$-Cocycle Superrigid action, and furthermore the corresponding orbit relations
are not compatible. More precisely:
\begin{thm}
\label{T:alg-not-quotient}
Let $H$ be a semi-simple Lie group, $\Delta<H$ an irreducible lattice, 
let $\Lambda<H$ be an unbounded countable subgroup acting by left translations
on $(H/\Delta,m_{H/\Delta})$. 

Then the algebraic ergodic p.m.p. action $\Lambda\acts (H/\Delta,m_{H/\Delta})$, is not
a quotient of any $\fintyp$-Cocycle Superrigid action.
In fact, if $\Gamma\acts(X,\mu)$ is a p.m.p. $\Delta$-Cocycle Superrigid action 
(e.g. a $\fintyp$-Cocycle Superrigid action),
of some countable group $\Gamma$,
then the orbit relations $\Rel{R}_{X,\Gamma}$ has no weak relation 
morphisms to  $\Rel{R}_{H/\Delta,\Lambda}$.  
\end{thm}
\begin{rems}
\label{R:Quotients-of-Bernoulli}
Let $\Gamma<G$ be a lattice in a simple a algebraic group $G$
of higher rank, and $\Gamma\acts (X,\mu)=(X_{0},\mu_{0})^{I}$ be an ergodic 
generalized Bernoulli action.
This action is $\fintyp$-Cocycle Superrigid by Popa's Cocycle Superrigidity
because such $\Gamma$ has property (T).
\begin{enumerate}
\item
It follows from Theorem~\ref{T:alg-not-quotient} above
that $\Gamma\acts (X,\mu)$ has no quotients of the form 
$\Gamma\acts (G/\Lambda,m_{G/\Lambda})$.
Thus in Theorem~\ref{T:ME+} case (i) applies. 
\item 
The proof of Theorem~\ref{T:alg-not-quotient}  is based on Popa's Cocycle Superrigidity.
However, one can use a direct argument to show that in this 
context $\Gamma\acts (X_{0},\mu_{0})^{I}$ has no non-trivial quotient
actions $\Gamma\acts (Y,\nu)$ where each element $g\in\Gamma$ has finite 
Kolmogorov-Sinai entropy: $h(Y,\nu,g)<\infty$ (see Section~\ref{S:Applications}).
Thus ME-Cocycle Superrigidity for Bernoulli actions of higher rank lattices
could have been proved before/independently of  Popa's Cocycle Superrigidity result.
\item
Using Theorem~\ref{T:superOE} the above statement can be strengthened to say that 
the orbit relation $\Rel{R}_{X,\Gamma}$ of a Bernoulli action of a higher rank lattice $\Gamma$
has no weak relation morphisms to orbit relations $\Rel{R}_{Y,\Lambda}$ of any essentially
free p.m.p. action $\Lambda\acts (Y,\nu)$ with element-wise finite entropy:
$h(Y,\nu,\ell)<\infty$ for all $\ell\in\Lambda$.
\end{enumerate}
\end{rems}


\bigskip

\subsection{Some more applications}
\label{SS:actions}
Theorem~\ref{T:PopaCSR} has some incidental applications to the 
structure of ergodic actions.
The following is immediate from the definitions:
\begin{prop}[{\bf Skew products}]
\label{P:skew}
Any skew-product of an $L$-Cocycle-Superrigid action $G\acts (X,\mu)$
by a p.m.p. action $L\acts (Z,\zeta)$ is isomorphic to a diagonal action on 
$(X\times Z,\mu\times\zeta)$ via some homomorphism $\ro:G\to L$.
\end{prop}
We denote by $\cpt\subset\fintyp$ the class of all separable compact groups.
\begin{prop}[{\bf Extension of weakly mixing malleable actions}]
\label{P:extensions}
Let $G\acts (X,\mu)$ be a $\cpt$-Cocycle Superrigid action
(e.g. a $\fintyp$-Cocycle Superrigid action), 
and $(X',\mu')\overto{p}(X,\mu)$ be a p.m.p. extension. 
If $G\acts (X',\mu')$ is a weakly mixing action then  $p$ is relatively weakly mixing extension.
\end{prop}
\begin{exmps}
\label{E:wmnotrelwm}
The following are examples of weakly mixing actions $\Gamma\acts (X,\mu)$
with a p.m.p. quotient $(X,\mu)\overto{p}(Y,\nu)$ which is not relatively weakly mixing;
in fact, here $p$ is a compact extension.
\begin{enumerate}
\item
Let $\Gamma\acts (X,\mu)=(X_{0},\mu_{0})^{I}$ be an ergodic 
generalized Bernoulli action of some group $\Gamma$. 
(Recall that such actions are weakly mixing).
Assuming that the base space $(X_{0},\mu_{0})$ is not an atomic probability space
with distinct weights, its automorphism group $\Aut(X_{0},\mu_{0})$ contains a non-trivial 
compact group $K$. The latter acts diagonally
on $(X,\mu)=(X_{0},\mu_{0})^{I}$ commuting with $\Gamma$. 
Thus $\Gamma\acts (X,\mu)$ has a quotient action $\Gamma\acts (Y,\nu)=(X,\mu)/K$,
and the natural projection $(X,\mu)\overto{p}(Y,\nu)$ is a compact extension
of $\Gamma$-actions.
(These actions exhibits very interesting rigidity phenomena
extensively studied by Popa in \cite{Popa:1cohom}, and by Popa and Vaes 
in \cite{PopaVaes:06}.)
\item
Let $\Gamma<G$ be a lattice in a non-compact center free simple Lie group $G$,
and $\Gamma\to K$ be a homomorphism into a compact group $K$ having  a 
dense image (e.g. $\Gamma=\SL_{n}(\ZZ)\to K=\SL_{n}(\ZZ_{p})$).
Then the $\Gamma$-action on $(K,m_{K})$ by translation is ergodic, and so is
the induced $G$-action on $X=(G\times K)/\Gamma\cong G/\Gamma\times K$.
By Howe-Moore theorem the action $G\acts X=G/\Gamma\times K$ is 
mixing, and hence weakly mixing. 
Yet it is a compact extension of the $G$-action on $Y=G/\Gamma$.
The restriction of these actions to any unbounded subgroup $\Gamma'<G$ 
(e.g. $\Gamma'=\Gamma$) 
inherits the above properties (another use of Howe-Moore).
\end{enumerate}
\end{exmps}

\begin{rems}
\label{R:BernoullibyK}
\begin{enumerate}
\item
In the above examples $\Gamma\acts (X,\mu)$ is a compact extension
of $\Gamma\acts (Y,\nu)$ defined via
some cocycle $\alpha:\Gamma\times Y\to K$ where $K$ is a compact group
with $X=Y\times K$. 
This cocycle $\alpha:\Gamma\times Y\to K$ cannot be untwisted in $K$,
because being weakly mixing the action $\Gamma\acts (X,\mu)$ 
admits no non-trivial compact quotients. 
\item
Taking $\Gamma\acts (X,\mu)$ to be $\fintyp$-Cocycle Superrigid
(e.g. Bernoulli action of a property (T) group $\Gamma$)
one gets a quotient action which fails to be  $\fintyp$-Cocycle Superrigid;
so relative weak mixing is a necessary assumption in Proposition \ref{P:CSRrelwm}.
\item
Note also that $\Gamma'\acts (G/\Gamma,m_{G/\Gamma})$ 
is not $\cpt$-Cocycle Superrigid, even though it is ME-Cocycle Superrigid
if $\Gamma$ and $\Gamma'$ are not commensurable.
\end{enumerate}
\end{rems}

Popa's Cocycle Superrigidity implies triviality of certain skew-product constructions
also in the context of equivalence relations:  
\begin{thm}
\label{T:semi-direct-rel}
Let $\Gamma\acts (X,\mu)$ be a $\fintyp$-Cocycle Superrigid action of a countable group $\Gamma$,
denote $\Rel{R}=\Rel{R}_{X,\Gamma}$, let $\Rel{S}$ some type ${\rm II}_{1}$ equivalence relation 
on some $(Y,\nu)$, and let $\Rel{Q}$ be an ergodic subrelation of  $\Rel{R}\times\Rel{S}$ on 
$(X\times Y,\mu\times\nu)$, such that the projection $X\times Y\to X$ is a relation
reduction of $\Rel{Q}$ to $\Rel{R}$. 

Then there exists a p.m.p. $\Gamma$-action $\Gamma\acts (Y,\nu)$ with 
$\Rel{R}_{Y,\Gamma}\subset\Rel{S}$, so that $\Rel{Q}$ differs from 
the orbit relation $\Rel{R}_{X\times Y,\Gamma}$ of the diagonal action 
$\Gamma$-action on $(X\times Y,\mu\times\nu)$ by an inner automorphism of
$\Rel{R}\times\Rel{S}$.
\end{thm}

\subsection*{Acknowledgments}
I am grateful to Sorin Popa for many enlightening conversations and 
for his helpful remarks on the early draft of these notes.
  
\section{Generalities and Preliminaries}

\subsection{$\Aut(X,\mu)$ as a topological group }
\label{SS:topologies}
Let $(X,\mu)$ be a Lebesgue measure space (finite or infinite). 
We denote by $\Aut(X,\mu)$ the group of all measure preserving 
measurable bijections between co-null sets of $X$, where two maps
which agree $\mu$-a.e. are identified. 
This group has two natural topologies: which we call \emph{uniform},
and \emph{weak} (also known as the vague or coarse topology).
The uniform topology can be defined by the metric
\begin{equation}
\label{e:strongtop}
	d^{(u)}(T,S)=\min\left(1,\,\mu\setdef{x\in X}{T(x)\neq S(x)}\right)
	\qquad(T,S\in\Aut(X,\mu)),
\end{equation}  
where the truncation by $1$ is relevant only to the infinite measure case.
Note that $d^{(u)}$ is a bi-invariant metric on $\Aut(X,\mu)$, which is easily 
seen to be complete. However, $\Aut(X,\mu)$ is not separable in this metric. 

We shall be almost exclusively using another, weaker, metrizable 
topology on $\Aut(X,\mu)$, with respect to which it becomes a Polish group,
i.e., complete second countable (equivalently separable) topological group. 
This topology is defined by a family of pseudo-metrics 
\begin{equation}
\label{e:weaktop}
	d^{(w)}_{E}(T,S)=\mu(TE \symdiff SE)
\end{equation}
where $E\subset X$ are measurable sets of finite measure.
One easily checks that $|d^{(w)}_{E}- d^{(w)}_{F}|\le 2\mu(E\symdiff F)$;
so the above topology is determined by countably many pseudo-metrics
$d^{(w)}_{E_{i}}$ where  $\{E_{i}\}_{i=1}^{\infty}$ are dense in the $\sigma$-algebra of $X$. 
Hence this topology is metrizable, the group is separable 
and is easily seen to be complete.

Both the uniform and the weak topologies are quite natural under the 
embedding $\pi:\Aut(X,\mu)\to U(L^{2}(X,\mu))$ in the unitary group of $L^{2}(X,\mu)$. 
Under this embedding the \emph{uniform} topology corresponds to the topology of the operator norm,
while the \emph{weak} topology corresponds to the strong/weak operator topology
(recall that the two coincide on the unitary group).
The fact that the image of $\Aut(X,\mu)$ is closed in $U(L^{2}(X,\mu))$ both in the 
weak/strong operator topology and in the norm topology follows from the completeness 
of the uniform and the weak topologies on $\Aut(X,\mu)$.

\subsection{Target groups -- classes $\fintyp\subset \biinv$}
\label{SS:Targetgroups}
In \cite{Popa:cocycle} Sorin Popa singles out a very interesting class of groups $\fintyp$.
By definition $\fintyp$  consists of all Polish groups which can be embedded
as closed subgroups of the unitary group of a finite von Neumann algebra. 
Groups of this class are well adapted to techniques involving 
unitary representations.  
Popa points out that the class $\fintyp$ is contained in the class $\biinv$ of all Polish groups
which admit a bi-invariant metric (see \cite[6.5]{Popa:cocycle}), and it seems to be unknown 
whether $\fintyp\subset\biinv$ is a proper inclusion.
We shall use the defining property of the class $\biinv$, which turns out to be very 
convenient for manipulations with abstract
measurable cocycles (see Section \ref{S:Cohom} below and forthcoming \cite{Furman:cohom}).
We shall consider the following subclasses as the main examples of groups in $\fintyp\subset\biinv$:
\begin{itemize}
\item 
	$\dscr$ -- the class of all discrete countable groups. $\dscr\subset \fintyp$ since
	a discrete group $\Lambda$ embeds in the unitary group of its 
	von Neumann algebra $L(\Lambda)$ which has type II$_{1}$. 
	The discrete metric on $\Lambda$ is clearly bi-invariant.
\item
	$\cpt$ -- the class of all second countable compact groups. 
	If $K$ is a compact separable group then its von Neumann algebra $L(K)$
	is of type II$_{1}$ (the trace is given by integration with respect to the Haar measure).
	To see directly that compact groups admit a bi-invariant metric, start from any 
	metric $d_{0}$ defining the 
	topology on $K$ and minimize, or average over $K\times K$-orbit to obtain an equivalent 
	bi-$K$-invariant metric:
	\[
		\qquad d_{inf}(g,h)\defq \inf_{K\times K} d_{k,k'}(g,h), \quad 
		d_{ave}(g,h)\defq\int_{K\times K} d_{k,k'}(g,h)\,dk\,dk'
	\]
	where $d_{k,k'}(g,h)=d_{0}(kgk', kh k')$, $(k,k'\in K)$.
\item
	$\Inn(\Rel{R})$	-- the inner automorphism group (a.k.a. as the \emph{full group}) of a 
	${\rm II}_{1}$ countable equivalence relation on a standard probability space $(X,\mu)$ 
	\[
		\Inn(\Rel{R})=\setdef{T\in\Aut(X,\mu)}{ (x,T(x))\in\Rel{R}\quad\text{a.e.}\quad x\in X}
	\] 
	Then $\Inn(\Rel{R})$ embeds in the unitary group of the  Murray - von Neumann algebra 
	$M(\Rel{R})$ associated to $\Rel{R}$; the type of the algebra is that of the the relation, i.e. 
	II$_{1}$ - a finite type. 
	
	More explicitly, $\Inn(\Rel{R})$ embeds as a closed subgroup  
	of $\Aut(\Rel{R},\tilde\mu)$, where $\tilde\mu$ denotes the infinite Lebesgue
	measure on $\Rel{R}$ obtained by lifting $\mu$ on $X$ via $p_{1}:\Rel{R}\to X$, 
	using the counting measure on the fibers. 
	Hence $\Inn(\Rel{R})$ inherits the weak topology from $\Aut(\Rel{R},\tilde\mu)$,
	coming from the weak/strong operator topology on the unitary group of
	$L^{2}(\Rel{R},\tilde\mu)$.
	The same topology can also be obtained as the restriction of the 
	uniform topology on $\Aut(X,\mu)$ to the subgroup $\Aut(\Rel{R})$ (see \ref{SS:topologies}).
	Since $d^{(u)}$ is a bi-invariant metric on $\Aut(X,\mu)$, in which $\Inn(\Rel{R})$
	is easily seen to be closed, this gives another explanation to the 
	inclusion $\Inn(\Rel{R})\in\biinv$.
\end{itemize}
Both classes $\fintyp\subset\biinv$ are closed under: passing to closed subgroups,
taking countable direct sums, taking direct integrals (i.e., passing from $G$ to $[0,1]^{G}$).

\subsection{Quotients and Skew products}
\label{SS:fibered}
A \emph{morphism} between two p.m.p. actions of $G$ on $(X,\mu)$ and $(Y,\nu)$ 
is a measurable map $p:X\to Y$ with $p_{*}\mu=\nu$, such that
for each $g\in G$: $p(g.x)=g.p(x)$ for $\mu$-a.e. $x\in X$.
We say that $Y$ is a \emph{quotient} of $X$, and/or that $X$ is an 
\emph{extension} of $Y$.
Given such a map one obtains a disintegration:
a measurable map $Y\ni y\mapsto \mu_{y}\in\Prob(X)$ such that
\[
	\mu=\int_{Y}\mu_{y}\,d\nu(y),\qquad \mu_{y}(p^{-1}(\{y\}))=1\quad
	\text{for}\quad \nu\text{-a.e. }y\in Y.
\]
Such a disintegration is unique (modulo null sets), and so it follows that for each $g\in G$:
$g_{*}\mu_{y}=\mu_{g.y}$ for $\nu$-a.e. $y\in Y$.

\medskip

Let $G\acts (Y,\nu)$ be a p.m.p. action, let $\sigma:G\times Y\to H$ 
be a measurable cocycle taking values in some Polish group $H$ which 
has a p.m.p. action $H\acts (Z,\zeta)$. 
Then one can construct the \emph{skew-product} $G$-action
\[
	(X,\mu)=(Y\times Z,\nu\times\zeta),\qquad g:(x,z)\mapsto (g.x, \sigma(g,x).z)
\]
For the skew-product action $G\acts (X,\mu)$ to be ergodic it is necessary, though
not sufficient, that both $G\acts (Y,\nu)$ and $V\acts (Z,\zeta)$ are ergodic. 

\medskip

We say that a morphism $(X,\mu)\to (Y,\nu)$ of $G$-actions is an 
\emph{isometric extension} if there exists a compact group $K$, a closed 
subgroup $K_{0}<K$, and a measurable cocycle $\sigma:G\times Y\to K$,
so that $G\acts (X,\mu)$ is measurably isomorphic to the skew product 
extension of $G\acts (Y,\nu)$
by $Z=K/K_{0}$ with the corresponding Haar measure. 

The $G$-action on $(Y,\nu)$ is a quotient of $G\acts (X,\mu)$; in fact, for ergodic
p.m.p. one can view any quotient $(X,\mu)\to (Y,\nu)$ of $G$-actions can be viewed
as  a skew-product defined by an action of $H< \Aut(Z,\zeta)$.
Here $(Z,\zeta)$ is a standard non-atomic Lebesgue space, unless $X\to Y$ is a finite
extension in which case we can take $Z=\{1,\dots,k\}$ and $H=S_{k}$.
The latter case is an example of an isometric extension of $Y$ using $K=S_{k}$ and 
$K_{0}\cong S_{k-1}$.



\subsection{Fibered Products}
\label{SS:skew}

Given two p.m.p. actions $G\acts (X_{i},\mu_{i})$, $i=1,2$,
which have a common quotient action $G\acts (Y,\nu)$, via $p_{i}:X_{i}\to Y$,
one can define the corresponding \emph{fibered product } space
\[
	X_{1}\times_{Y}X_{2}=\setdef{(x_{1},x_{2})\in X_{1}\times X_{2}}{p_{1}(x_{1})=p_{2}(x_{2})}
\]
equipped with a fibered product probability measure
\[
	\mu_{1}\times_{\nu}\mu_{2}=\int_{Y}\mu_{1,y}\times\mu_{2,y}\,d\nu(y)
\]
where $\mu_{i}=\int \mu_{i,y}\,d\nu(y)$ are the corresponding disintegrations.
The measure $\mu_{1}\times_{Y}\mu_{2}$ is invariant under the diagonal $G$-action 
$g:(x_{1},x_{2})\mapsto (g.x_{1},g.x_{2})$ on $X_{1}\times_{Y}X_{2}$.
The definition of fibered products has an obvious extension to more than two factors.
We shall denote by $(X^{n}_{Y},\mu^{n}_{\nu})$ the fibered product
of $1\le n\le\infty$ copies of a given $G$-action $G\acts (X,\mu)$ with the quotient $G$-action
on $(Y,\nu)$. 
The usual product (which can be viewed as the fibered product over the trivial $Y$) 
will be denoted by $(X^{n},\mu^{n})$.
\subsection{Relative weak mixing}
\label{SS:relwm}
Recall that a p.m.p. action $G\acts (X,\mu)$ is \emph{weakly mixing} if the diagonal $G$-action 
on $(X\times X,\mu\times\mu)$ is ergodic. 
If $G\acts (X,\mu)$ is weakly mixing then the diagonal $G$-action on $(X\times X',\mu\times\mu')$
is ergodic for any ergodic p.m.p. action $G\acts (X',\mu')$, and it follows that 
the  diagonal $G$-action on the product $(X^{n},\mu^{n})$ is ergodic, and weakly mixing.
 
We shall be interested in the relative version of the notion of weak mixing, introduced
by Furstenberg in \cite{Furstenberg:Szemeredi:77}, 
and  Zimmer in \cite{Zimmer:IllJ1:76}, \cite{Zimmer:IllJ2:76}. 
Glasner's \cite[Ch 9,10]{Glasner:Joinings:03} is a good recent reference.
\begin{thm}[Relative weak mixing]
\label{T:relwmcond}
Let $(X,\mu)\to (Y,\nu)$ is a morphism of ergodic $G$-actions. The following are equivalent:
\begin{enumerate}
\item
	The morphism is relatively weakly mixing, i.e., $G\acts (X^{2}_{Y},\mu^{2}_{\nu})$ is ergodic.
\item
	For any ergodic p.m.p. action $G\acts (X',\mu')$ which has $G\acts (Y,\nu)$
	as a quotient, the $G$-action on  $(X\times_{Y}X',\mu\times_{\nu}\mu')$ is ergodic.
\item
	The diagonal $G$-action on $(X^{n}_{Y},\mu^{n}_{\nu})$ is ergodic for all $n\ge 2$.
\item
	There does not exist intermediate quotients $(X,\mu)\to(Y',\nu')\to (Y,\nu)$ of the $G$-actions, 
	where $(Y',\nu')\to (Y,\nu)$ is a non-trivial isometric extension.
\end{enumerate}
\end{thm}
Of course $G\acts (X,\mu)$ is weakly mixing iff it is weakly mixing relative to the trivial 
action on a point.
\begin{warn}
If $(X,\mu)\to (Y,\nu)\to (Z,\zeta)$ are morphisms of ergodic $G$-actions
and $(X,\mu)\to (Z,\zeta)$ is relatively weakly mixing, then
$(Y,\nu)\to (Z,\zeta)$ is also relatively weakly mixing, 
but $(X,\mu)\to (Y,\nu)$ need not be such.
In particular, a weakly mixing action $G\acts (X,\mu)$
can fail to be weakly mixing relative to some quotients
(recall Examples \ref{E:wmnotrelwm}).
\end{warn}

\subsection{Relative Property (T)}
\label{SS:relT}
Let $\pi:G\to U(\sH_\pi)$ be a unitary representation of some topological group
(unitary representations are always assumed to be continuous with respect to 
the weak/strong topology on the unitary group).
Given a compact subset $K\subset G$ and $\epsilon>0$ let
\[
	V_{K,\epsilon}=\left\{v\in\sH_{\pi}\left| 
	\|v\|=1,\ \max_{g\in K}\|\pi(g)v-v\|<\epsilon\right.\right\}
\]
Vectors in $V(K,\epsilon)$ are called $(K,\epsilon)$-almost $\pi(G)$-invariant.
A unitary $G$-representation $\pi$ is said to \emph{almost have invariant vectors}
if $V_{K,\epsilon}$ is non-empty for any compact $K\subset G$ and $\epsilon>0$. 

A \lcsc group $G$ has Kazhdan's property (T) if every unitary $G$-representation
which \emph{almost has invariant vectors}, actually has non-trivial invariant vectors.
A closed subgroup $H<G$ is said to have relative property (T) in $G$, if every 
unitary $G$-representation $\pi$ which almost has $G$-invariant vectors,
has $H$-invariant vectors.

The fact that a sequence of unitary representations can be organized into a single one, by taking
their direct sum, can be used to show that the quantifiers in the definition of property (T)
can be interchanged to say that: $G$ has property (T) (resp. $H<G$ has relative (T))
if there exist compact $K\subset G$ and $\epsilon>0$, so that any unitary 
$G$-representation with non-empty $V_{K,\epsilon}$ has non-trivial $G$-invariant vectors
(resp. $H$-invariant vectors). 
Furthermore, the invariant vectors can be found close to the almost invariant ones:
\begin{thm}[Jolissaint \cite{Jolissaint:relT:05}]
\label{T:alm-inv}
If $H<G $ has relative property (T), then for any $\delta>0$ there exists a 
compact subset $K\subset G$ and $\epsilon>0$ such that: 
any unitary $G$-representation $\pi$ with non-empty $V_{K,\epsilon}$
has an $H$-invariant unit vector $u$ with $\dist(u,V_{K,\epsilon})<\delta$.
\end{thm}
For property (T) groups (i.e. the case where $H=G$) or, more generally,
in the case where $H$ is normal in $G$, the above result is easy 
(consider the $G$-representation on the orthogonal compliment of the space
of $H$-invariant vectors). 
For the general case of an arbitrary subgroup $H<G$ with relative property (T),
a different argument is needed (see \cite{Jolissaint:relT:05}, or the appendix in 
\cite{PetersonPopa:relT:05}).

Many property (T) groups are known (see \cite{dHarpeValette:T:89}, \cite{Zuk:GAFA:03}),
but there are many many more groups with normal subgroups with relative property (T).
Basic example of this type is $H=\ZZ^2$ in $G=\SL_2(\ZZ)\ltimes\ZZ^2$.
Many arithmetic examples are constructed in Alain Valette's \cite{Valette:relT:05}.
A very general construction of relative property (T) inclusions of an Abelian
$H$ into semi-direct products $G=\Gamma\ltimes H$ were recently obtained by 
Talia Fern\'os in \cite{Fernos:relT}.

\subsection{Gaussian Actions}
\label{SS:Gaussian}
Let us briefly describe this basic construction of p.m.p. actions stemming from unitary (or rather orthogonal)
representations. For references see \cite{ConnesWeiss:T:80}, \cite{Glasner:Joinings:03},
\cite{Cherixetal:01} and references therein.

The standard Gaussian distribution on $\RR$ has density $(2\pi)^{-1/2}e^{-x^{2}/2}\,dx$;
the direct product of $k$ copies of this measure, denoted $\mu_{\RR^{k}}$,
is the  \emph{standard Gaussian} measure on $\RR^{k}$; it has density
$(2\pi)^{-k/2}e^{-\|x\|^{2}/2}\,dx$ on $\RR^{k}$. 
The standard measure $\mu_{\RR^{k}}$ is invariant under orthogonal transformations, 
and orthogonal projections $\RR^{k}\to\RR^{l}$ map $\mu_{\RR^{k}}$ to $\mu_{\RR^{l}}$.
A general (centered) Gaussian measure on a finite dimensional vector space $\RR^{n}$
is a pushforward $A_{*}\mu_{\RR^{k}}$ of the standard Gaussian measure $\mu_{\RR^{k}}$
by a linear map $A:\RR^{k}\to\RR^{n}$ (one can always take $k\le n$ and $A=UD$
where $D$ is diagonal and $U$ is orthogonal).
For $\mu=A_{*}\mu_{\RR^{k}}$ the coordinate projections 
have correlation matrix $C=AA^{T}$; in fact a general (centered) Gaussian measure on $\RR^{n}$
is determined by its correlation matrix $C$ and therefore will be denoted $\mu_{C}$.
If $C$ is invertible, we say that $\mu_{C}$ is non-degenerate; in this case 
$d\mu_{C}(y)=(2\pi)^{-n/2}\det(C)^{-1}e^{-\ip{C^{-1}y}{y}/2}\,dy$.

Let $I$ be a countable set, $\{v_{i}\}_{i\in I}$ a sequence of unit vectors in some real Hilbert space $\sH$,
and $C:I\times I\to[-1,1]$ be given by $C_{i,j}=\ip{v_{i}}{v_{j}}$.
Then for any finite subset $J=\{i_{1},\dots,i_{n}\}\subset I$ the $n\times n$ matrix
$C_{J}\defq C|_{J\times J}$ is a correlation matrix, and defines a centered Gaussian
probability measure $\mu_{C_{J}}$ on $\RR^{J}$.  
It is easy to see that for any two finite subsets $K\subset J$ of $I$ the projection of $\mu_{C_{J}}$ 
under $\RR^{J}\to \RR^{K}$ is $\mu_{C_{K}}$. 
Hence by Kolmogorov's extension theorem there exists
a unique probability measure $\mu_{C}$ on $\RR^{I}$ whose projections on $\RR^{J}$ are $\mu_{C_{J}}$
for any finite subset $J\subset I$.

Now let $\pi:G\to O(\sH)$ be a  linear orthogonal representation of a topological group $G$,
and assume that $G$ leaves invariant the set $\{v_{i}\}_{i\in I}$. Then for any $g\in G$:
\[
	\ip{\pi(g)v_{i}}{\pi(g)v_{j}}=\ip{v_{i}}{v_{j}},\qquad (i,j\in I)
\] 
and therefore the map $g:\RR^{J}\to\RR^{gJ}$ maps $\mu_{C_{J}}$ onto $\mu_{C_{gJ}}$.
This defines a measure preserving $G$-action on $(\RR^{I},\mu_{C})$, which is called the 
\emph{Gaussian process/action} associated to the unitary $G$-action on the set $V=\{v_{i}\}_{i\in I}$. 
Assuming $V$ spans a dense subspace of $\sH$, we get an isometric embedding $\sH\to L^{2}(\RR^{I},\mu_{C})$
sending $v_{i}$ to the coordinate projection $f_{i}:\RR^{I}\to\RR$.
In particular $\pi$ becomes a subrepresentation of the $G$-representation on $L^{2}(\RR^{I},\mu_{C})$.
\footnote{Often Gaussian actions are described as the action on a probability space on which
a generating family $\{X_{i}\}_{i\in I}$ of standard Gaussian variables is defined with correlations 
${\rm E}(X_{i}X_{j})=\ip{v_{i}}{v_{j}}$; the $G$-action corresponds to the permutation of $X_{i}$s.
In our description the underlying space is $(\RR^{I},\mu_{C})$ and $X_{i}$s are coordinate projections.}

\begin{rem}
\label{R:Bernoulli-is-Gaussian}
Observe that Bernoulli actions with a non-atomic base space is a particular case
of a Gaussian process, corresponding to $C$ being the identity matrix, i.e., $\{v_{i}\}_{i\in I}$
being orthonormal, the situation arising in the $G$-representation on $\sH=\ell^{2}(I)$.
\end{rem}

Gaussian actions can be defined not only for representations of discrete groups, but for continuous
orthogonal repesentations $\pi:G\to O(\sH)$ of an arbitrary \lcsc group $G$. 
The fastest explaination is the following: given $\pi:G\to O(\sH)$ of such a group $G$ 
choose a dense countable subgroup $\Lambda<G$
and apply the above construction to a $\pi(\Lambda)$-invariant spanning set 
$\{v_{i}\}_{i\in I}\subset \sH$.
Note that continuity of the representation $\pi$ implies continuity of the $\Lambda$-action on 
$(\RR^{I},\mu_{C})$, i.e., if $\lambda_{n}\to 1$ in $G$ then 
\begin{equation}
\label{e:cont}
	\mu_{C}(\lambda_{n}E\triangle E)\to 0,\qquad(E\in\mathcal{B}(\RR^{I})).
\end{equation} 
Indeed, it suffices to check this for sets $E$ which appear as a preimage of a
Borel subset of $\RR^{J}$ of the projection $\RR^{I}\to \RR^{J}$, where $J\subset I$ is finite.
For such sets (\ref{e:cont}) follows from the fact that 
$\lim_{n\to \infty}\max_{i\in J}\|\pi(g_{n})v_{j}-v_{j}\|= 0$.
Hence, the $\Lambda$-action on $(X,\mu)\defq(\RR^{I},\mu_{C})$ extends to 
$G\to \Aut(X,\mu)$ by continuity.  
Since $G$ is \lcsc this homomorphism is realized as a measurable action 
$G\acts (X,\mu)$.

\section{Cohomology of Measurable Cocycles}
\label{S:Cohom}
This section contains some results about measurable cocycles $G\times X\to L$
which do not depend on property (T) of $G$ or on malleability of the action $G\acts (X,\mu)$ 
(results related to the latter notions are discussed in Section~\ref{S:rigidity-deformation}).
Hence the nature of the acting \lcsc group $G$ is immaterial. 
However, it will be important to impose certain conditions on the target group $L$.
Here we will assume that $L\in \biinv$, namely that it admits a bi-invariant metric $d$
with respect to which it is a complete separable group. 
The importance of this condition and of the weak mixing assumption 
will become clear in the following Lemma~\ref{L:basic}, which is a basis for the 
main statements \ref{L:target} -- \ref{L:normal} of this section.
These statements parallel Lemma 2.11, Proposition 3.5, Theorem 3.1, 
and Lemma 3.6 of Popa's \cite{Popa:cocycle}, where they are proved in the
context of Operator Algebras. 
Here we work in the framework of the class $\biinv\supset\fintyp$ and use
only elementary Ergodic theoretic arguments
We remark that similar results hold for some other classes of target groups
and, sometimes without weak mixing assumption.
This is discussed in the forthcoming paper \cite{Furman:cohom}. 

We start by an elementary observation:
\begin{lem}[{\bf Separability argument}]
\label{L:separability}
Let $(X,\mu)$ be a probability space, $(M,d)$ a separable metric space,
$\Phi:X\to M$ a Borel map such that $d(\Phi(x_{1}),\Phi(x_{2})=d_{0}$ for 
$\mu\times\mu$-a.e. $(x_{1},x_{2})\in X\times X$.
Then $d_{0}=0$ and $\Phi_{*}\mu=\delta_{m_{0}}$ for some $m_{0}\in M$. 
\end{lem}
\begin{proof}
Assume $d_{0}>0$. By separability, $M$ can be covered by countably many open balls 
$M=\bigcup_{1}^{\infty} B_{i}$ with $\diam(B_{i})=\sup\{d(p,q) \mid p,q\in B_{i}\}<d_{0}$. 
Note that $E_{i}=\Phi^{-1}(B_{i})$ has $\mu(E_{i})=0$, but this is impossible 
because $1=\sum_{i=1}^{\infty} \mu_{y}(E_{i})$. 
Hence $d_{0}=0$ and therefore $\Phi_{*}\mu$ is a Dirac measure $\delta_{m_{0}}$ at
some $m_{0}\in M$. 
\end{proof}

The following should be compared to Popa's original \cite[Lemma 2.11]{Popa:cocycle}.
\begin{lem}[{\bf Basic}]
\label{L:basic}
Let $G\acts (X,\mu)$ be an ergodic action, $(X,\mu)\overto{p} (Y,\nu)$ 
a $G$-equivariant measurable quotient which is relatively weakly mixing.
Let $L$ be a group in class $\biinv$, and let 
$\alpha,\beta:G\times Y\to L$ be two measurable cocycles.
Let $\Phi:X\to L$ be a measurable function, so that for each $g\in G$
for $\mu$-a.e. $x\in X$:
\[
	\alpha(g,p(x))=\Phi(g.x) \beta(g,p(x)) \Phi(x)^{-1}.
\]
Then $\Phi$ descends to $Y$: 
there exists a measurable $\phi:Y\to L$ so that $\Phi=\phi\circ p$ a.e. on $X$, 
and  
\[
	\alpha(g,y)=\phi(g.y) \beta(g,y) \phi(y)^{-1}.
\]
\end{lem}
This Lemma describes the \emph{injectivity} of the pull-back map 
\[
	H^{1}(G\times X,L)\overfrom{p^{*}} H^{1}(G\times Y,L)
\] 
between the sets of equivalence classes of measurable cocycles, 
corresponding to weakly mixing morphisms $X\overto{p} Y$ 
and coefficients $L$ from $\biinv$. 
This implies Proposition~\ref{P:CSRrelwm} in the introduction.
\begin{rem}
The assumption of relative weak mixing is essential
(recall Remark~\ref{R:BernoullibyK}). 
For ergodic extensions $(X,\mu)\overto{p}(Y,\nu)$ which are not necessarily 
relatively weakly mixing one can prove that if 
cocycles $\alpha,\beta:G\times Y\to L$ (with $L\in\biinv$) are cohomologous 
when lifted to $G\times X\to L$, then they are already cohomologous when lifted
to an intermediate isometric extension $(Y',\nu')\to(Y,\nu)$
(a finite extension, if $L$ is countable).
This is discussed in more detail in \cite{Furman:cohom}.
\end{rem}
\begin{proof}
Consider the diagonal $G$-action on the fibered product 
$(X\times_{Y}X,\mu\times_{\nu}\mu)$. 
It is ergodic by the assumption of relative weak mixing.
The function $f:X\times_{Y}X\to[0,\infty)$ defined by 
\[
	f(x_{1},x_{2})=d(\Phi(x_{1}),\Phi(x_{2}))
\] 
is measurable and is $G$-invariant. 
Indeed for any $g\in G$ for $\mu\times_{\nu}\mu$-a.e.
$(x_{1},x_{2})$ we have, denoting by $y=p(x_{1})=p(x_{2})$:
\begin{eqnarray*}
	f(g.x_{1},g.x_{2})&=&d(\Phi(g.x_{1}),\Phi(g.x_{2}))\\
	&=&d(\alpha(g,y) \Phi(x_{1}) \beta(g,y)^{-1},\alpha(g,y)\Phi(x_{2})\beta(g,y)^{-1})\\
	&=&d(\Phi(x_{1}),\Phi(x_{2}))=f(x_{1},x_{2}).
\end{eqnarray*} 
Ergodicity of the $G$-action on the fibered product (i.e. the assumption of relative weak mixing) 
implies that $f$ is essentially a constant $d_{0}$, i.e. $f(x_{1},x_{2})=d_{0}$ for 
$\mu\times_{\nu}\mu$-a.e. $(x_{1},x_{2})$.
For $\nu$-a.e. $y\in Y$ we apply Lemma~\ref{L:separability} to $(X,\mu_{y})$ to deduce
that $\Phi_{*}\mu_{y}=\delta_{\phi(y)}$ for some $\phi:Y\to L$.
Measurability of $\phi$ follows from that of $\Phi$.
\end{proof}

\bigskip

The following is a relative version of Popa's \cite[Theorem 3.1]{Popa:cocycle}.
Our proof is based on the Basic Lemma above. 
\begin{thm}[{\bf Untwisting cocycles}]
\label{T:untwisting}
Let $G\acts (X,\mu)$ be p.m.p. action which is weakly mixing relative to a quotient 
action $G\acts (Y,\nu)$;
let $L$ be a group from class $\biinv$, and let $\alpha,\beta :G\times X\to L$ be two 
measurable cocycle.
Assume that there exists a map $F:X\times_{Y} X\to L$ so that
for every $g\in G$: for $\mu^{2}_{\nu}$-a.e. $(x_{1},x_{2})\in X^{2}_{Y}$ 
\begin{equation}
\label{e:abF}
	\alpha(g,x_{1})=F(g.x_{1}, g.x_{2})\,\beta(g,x_{2})\,F(x_{1},x_{2})^{-1}.
\end{equation}
Then there exists a measurable cocycle $\ro:G\times Y\to L$ and measurable maps 
$\phi,\psi:X\to L$ so that
\begin{eqnarray*}
	\alpha(g,x)&=&\phi(g.x)\,\ro(g,p(x))\,\phi(x)^{-1},\\
	\qquad \beta(g,x)&=&\psi(g.x)\,\ro(g,p(x))\,\psi(x)^{-1}
\end{eqnarray*}
for all $g\in G$ and $\mu$-a.e. $x\in X$.
\end{thm}
Of course the ``absolute case'' corresponds to $G\acts (Y,\nu)$ being the trivial 
action on a point; in this case $\ro:G\to L$ is a homomorphism. 
The reader, interested in this case only, should read the following short proof
consistently replacing the phrase ``for $\nu$-a.e. $y\in Y$ for $\mu_{y}$-a.e. ... '' 
by just ``for $\mu$-a.e. ....''.
\begin{proof}
By the assumption for $\nu$-a.e. $y\in Y$ the relation (\ref{e:abF}) holds
for $\mu_{y}\times\mu_{y}$-a.e. $(x_{1},x_{2})$. By Fubini for $\nu$-a.e. $y\in Y$ for
$\mu_{y}\times\mu_{y}\times\mu_{y}$-a.e. $(x_{1},x_{2},x_{3})$ we have both relations:
\begin{eqnarray*}
		\alpha(g,x_{1})&=&F(g.x_{1},g.x_{2})\,\beta(g,x_{2})\,F(x_{1},x_{2})^{-1}
		\qquad\text{and}\\
		\alpha(g,x_{3})&=&F(g.x_{3},g.x_{2})\,\beta(g,x_{2})\,F(x_{3},x_{2})^{-1}.
\end{eqnarray*}
Substituting the  first in the second, we obtain the following identities which
hold $\mu^{3}_{\nu}$-almost everywhere on $X^{3}_{Y}$:
\begin{eqnarray*}
	\alpha(g,x_{3})&=&F(g.x_{3},g.x_{2})\,\beta(g,x_{2})\,F(x_{3},x_{2})^{-1}\\
			&=&F(g.x_{3},g.x_{2})\,F(g.x_{1},g.x_{2})^{-1}\,\alpha(g,x_{1})\,
				F(x_{1},x_{2})\,F(x_{3},x_{2})^{-1}.
\end{eqnarray*}
Setting $\Phi(x_{1},x_{2},x_{3})\defq F(x_{1},x_{2})\,F(x_{3},x_{2})^{-1}$ the above takes the form
\[
	\alpha(g,x_{1})=\Phi(g.x_{1},g.x_{2},g.x_{3})\,\alpha(g,x_{3})\,\Phi(x_{1},x_{2},x_{3})^{-1}
\]
or, equivalently:
\[
	\Phi(g.x_{1},g.x_{2},g.x_{3})=\alpha(g,x_{1})\,\Phi(x_{1},x_{2},x_{3})\,\alpha(g,x_{3})^{-1}.
\]
Next, observe that the morphism of the diagonal $G$-actions 
\[
	q:(X^{3}_{Y},\mu^{3}_{\nu})\to (X^{2}_{Y},\mu^{2}_{\nu}),\qquad 
	q(x_{1},x_{2},x_{3})=(x_{1},x_{3})
\]
is relatively weakly mixing. Indeed, weak mixing of this morphism is equivalent
to the ergodicity of the diagonal $G$-action on the fibered
product 
\[
	X^{3}_{Y}\times_{(X^{2}_{Y})} X^{3}_{Y}\cong X^{4}_{Y}
\] 
Ergodicity of the $G$-action on the latter
follows from relative weak mixing (Theorem~\ref{T:relwmcond} (1)$\Longrightarrow$ (3)).

Hence, by the Basic Lemma~\ref{L:basic},  $\Phi(x_{1},x_{2},x_{3})=f(x_{1},x_{3})$
for some measurable map $f:X^{2}_{Y}\to L$. Therefore 
\[
	F(x_{1},x_{2})\,F(x_{3},x_{2})^{-1} = \Phi(x_{1},x_{2},x_{3})=f(x_{1},x_{3})
\]
meaning that $\mu^{3}_{\nu}$-a.e. on $X^{3}_{Y}$: 
\[
	F(x_{1},x_{2})=f(x_{1},x_{3})F(x_{3},x_{2}).
\]
Using Fubini one can choose a measurable section $s:Y\to X$ so that
for $\nu$-a.e. $y\in Y$ for $\mu_{y}\times\mu_{y}$-a.e. $(x_{1},x_{2})$:
\[
	F(x_{1},x_{2})=f(x_{1},s(y)) F(s(y),x_{2}).
\]
Defining $\phi,\psi:X\to L$ by
\[
	\phi(x)\defq f(x,s\circ p(x)),\qquad \psi(x)\defq F(s\circ p(x),x)^{-1}
\]
we get that $F(x_{1},x_{2})$ splits as $F(x_{1},x_{2})=\phi(x_{1})\psi(x_{2})^{-1}$.
This allows to rewrite (\ref{e:abF}) as $\mu^{2}_{\nu}$-a.e. identity:
\[
	\phi(g.x_{1})^{-1}\,\alpha(g,x_{1})\,\phi(x_{1})=\psi(g.x_{2})^{-1}\,\beta(g,x_{2})\,\psi(x_{2}).
\]
Hence for each $g\in G$:  for $\nu$-a.e. $y\in Y$ the above holds
for $\mu_{y}\times\mu_{y}$-a.e. pair $(x_{1},x_{2})$. 
Therefore the left and the right hand sides are each $\mu_{y}$-a.e. constant, which we denote
by $\ro(g,y)$. 
This measurable function $\ro:G\times Y\to L$ inherits the cocycle
equation from $\alpha$ and/or $\beta$.
\end{proof}

\bigskip

In the proof of the Popa Cocycle Superrigidity theorem in the case of group
inclusion $H<G$ with relative property (T), one first shows that the restriction of
a cocycle $\alpha:G\times X\to L$ to $H\times X\to L$ can be ``untwisted'' 
so to descend to a cocycle $\ro:H\times Y\to L$ 
(or to become a homomorphism $\ro:H\to L$ in the ``absolute'' case).
One then needs to show that this is enough to untwist the whole cocycle $\alpha:G\times X\to L$.
The next Lemma, parallel to Popa's \cite[Lemma 3.6]{Popa:cocycle} 
(going back to  \cite[Lemma 5.7]{Popa:noncommB:06}), shows how the 
``untwisting'' on $H$ extends to a larger subgroup $\hat{H}<G$.
\begin{lem}[{\bf Promoting homomorphisms}]
\label{L:normal}
Let $G\acts (X,\mu)$ be a p.m.p. action with a quotient $(Y,\nu)$,
$H<G$ be subgroups, and $\alpha:G\times X\to L$, $\ro:H\times Y\to L$
be measurable cocycles taking values in a group $L$ from class $\biinv$.
Suppose that the restriction of $\alpha$ to $H\times X$ descends to $\ro$, namely:
\[
	\alpha(h,x)=\ro(h,p(x))\qquad (h\in H).
\]
and that $H\acts (X,\mu)$ is relatively weakly mixing relative to $(Y,\nu)$.

Then, denoting $\hat{H}=N_{G}(H)$, the cocycle $\ro$ extends to 
$\hat\ro:\hat{H}\times Y\to L$ and the restriction of $\alpha$ to $\hat{H}\times X$ 
descends to $\hat\ro$:
\[
	\alpha(g,x)=\hat\ro(g,p(x))\qquad (g\in \hat{H}).
\]
Moreover, if there exists a subgroup $K<H$ so that $K\acts (X,\mu)$ is weakly
mixing relative to $(Y,\nu)$, then the same conclusion applies to the group
$\hat{H}$ generated by the set $G_{K,H}\defq\setdef{g\in G}{gKg^{-1}<H}\supseteq N_{G}(H)$.
\end{lem}
\begin{proof}
Let $K<H$ (possibly $K=H$) be a subgroup so that $K\acts (X,\mu)$ is 
weakly mixing relative to $(Y,\nu)$.
For $g\in G$ we denote by $k\mapsto k^{g}\defq gkg^{-1}$ the conjugation  
isomorphism $K\overto{\cong}K^{g}=\setdef{k^{g}}{k\in K}$. 

Fix some $g\in G_{K,H}$. Computing $\alpha(gk,x)=\alpha(k^{g}g,x)$ 
using the cocycle identity gives (using that $K^{g}<H$):
\begin{eqnarray*}
	\alpha(gk,x) &=& \alpha(g,k.x)\,\alpha(k,x)=\alpha(g,k.x)\,\ro(k,y)\\
	\alpha(k^{g}g,x)&=&\alpha(k^{g},g.x)\,\alpha(g,x)=\ro(k^{g},g.y)\,\alpha(g,x)
\end{eqnarray*}
where $y=p(x)$. This gives a $\mu$-a.e. identity
\[
	\alpha(g,k.x)=\ro(k^{g},g.y)\,\alpha(g,x)\,\ro(k,y)^{-1}.
\]
Observe that the map $\ro^{g}:K\times Y\to L$ defined by $\ro^{g}(k,y)\defq \ro(k^{g},g.y)$ 
is a cocycle.
Therefore using Basic Lemma \ref{L:basic} the map $\Phi(x)=\alpha(g,x)$, satisfying
$
	\Phi(k.x)=\ro^{g}(k,y)\,\Phi(x)\,\ro(k,y)^{-1},
$
descends to a function $Y\to L$. The set 
\[	
	G_{1}=\setdef{g\in G}{\alpha(g,x)\ \text{depends only on } y=p(x)}
\] 
forms a subgroup of $G$. 
The above arguments shows that $G_{1}\supseteq G_{K,H}\supseteq N_{G}(H)$.
Hence $G_{1}>\hat{H}$ and $\hat{\ro}(g,y)$ denotes the value of $\alpha(g,x)$, $y=p(x)$, 
for $g\in G_{1}>\hat{H}$.
\end{proof}

The following Lemma is parallel to Popa's \cite[Proposition 3.5]{Popa:cocycle}.
We shall refer to this Lemma only in Remark \ref{R:general-rigidity}.
\begin{lem}[{\bf Tightness for untwisting}]
\label{L:target}
Let $G\acts (X,\mu)$ be an ergodic action, $(X,\mu)\overto{p} (Y,\nu)$ 
a $G$-equivariant measurable quotient which is relatively weakly mixing,
$L$ a group in class $\biinv$, $M<L$ a closed subgroup,
$\alpha:G\times X\to M$ and $\ro:G\times Y\to M$ be measurable cocycles.
Suppose that there exists a measurable map $\Phi:X\to L$ so that for all $g\in G$ for $\mu$-a.e. $x$:
\begin{equation}
\label{e:Lcohom}
	\ro(g,p(x))=\Phi(g.x)\,\alpha(g,x)\,\Phi(x)^{-1}.
\end{equation}
Then there exists a measurable map $\Phi':X\to M$ and a cocycle
$\ro':G\times Y\to M$ so that for all $g\in G$ for $\mu$-a.e. $x$:
\begin{equation}
\label{e:Mcohom}
	\ro'(g,p(x))=\Phi'(g.x)\,\alpha(g,x)\,\Phi'(x)^{-1}.
\end{equation}
In the particular case of trivial $Y$, i.e., if $G\acts (X,\mu)$ is weakly mixing 
and the cocycle $\alpha:G\times X\to M$ can be untwisted to a homomorphism 
$\ro:G\to M<L$ in $L$, then $\alpha$ 
can be untwisted to a homomorphism $\ro':G\to M$ within $M$.
\end{lem}
\begin{proof}
The assumption that $L\in\biinv$ is used to define the following
$L$-invariant metric $\bar{d}$ on $L/M$:
\[
	\bar{d}(\ell_{1}M,\ell_{2}M)=\dist(\ell_{2}^{-1}\ell_{1},M)
	=\inf_{m\in M} d(\ell_{2}^{-1}\ell_{1},m).
\]
Let $f:X\times_{Y} X\to [0,\infty)$ be given by
$
	f(x_{1},x_{2})=\bar{d}(\Phi(x_{1})M,\Phi(x_{2})M).
$
This measurable function is invariant for the diagonal $G$-action on 
$X\times_{Y}X$ because 
\begin{eqnarray*}
	f(g.x_{1},g.x_{2})&=&\bar{d}(\Phi(g.x_{1})M,\Phi(g.x_{2})M)\\
	&=&\bar{d}(\ro(g,y) \Phi(x_{1})\alpha(g,x)^{-1}M,\ro(g,y)\Phi(x_{2})\alpha(g,x)^{-1}M)\\
	&=&\bar{d}(\ro(g,y) \Phi(x_{1})M,\ro(g,y)\Phi(x_{2})M)\\
	&=&\bar{d}(\Phi(x_{1})M,\Phi(x_{2})M)=f(x_{1},x_{2})
\end{eqnarray*} 
where $y=p(x)\in Y$.
Ergodicity of the $G$-action on the fibered product (i.e. the assumption of relative weak mixing) 
implies that $f(x_{1},x_{2})$ is $\mu\times_{\nu}\mu$-a.e. constant $d_{0}$.
Using Lemma~\ref{L:separability} we deduce $d_{0}=0$, and so for 
some measurable $\phi:Y\to L$ we have
\[
	\Phi(x)M=\phi(\pi(x))M,
\]
which implies that $\Phi'(x)=\phi(\pi(x))^{-1}\Phi(x)$ takes values in $M$.
Let $\ro'(g,y)=\phi(g.y)^{-1}\ro(g,y)\phi(y)$. This is a cocycle on $G\times Y$ satisfying
(\ref{e:Mcohom}). All the terms n the right hand side of this equation take
values in $M$ and so does $\ro':G\times Y\to M$ as claimed.
\end{proof}

\section{Rigidity vs. Deformation} 
\label{S:rigidity-deformation}
\subsection{The topology on $Z^{1}(G\times X, L)$}
\label{D:top-on-coc}
Let $G\acts (X,\mu)$ be a fixed p.m.p. action of a \lcsc group, let $L$ be some Polish group
(which will be assumed to be in class $\biinv$). 
Consider the space $Z^{1}(G\times X,L)$ of all measurable cocycles $G\times X\to L$
(identified modulo $\mu$-null sets)
with the topology of convergence in measure, which can be defined as follows. 
Fix a left-invariant (below bi-invariant) metric $d$ on $L$.
Given a compact subset $K\subset G$ and $\epsilon>0$ let $V_{K,\epsilon}$ 
denote the set of pairs $(\alpha,\beta)$ of cocycles $\alpha,\beta:G\times X\to L$ such that
\[
	\mu\setdef{x\in X}{d(\alpha(g,x),\beta(g,x))<\epsilon}>1-\epsilon
	\qquad(g\in K).
\]
Sets $V_{K,\epsilon}$ form the base of topology on the cocycles. This topology is
complete and is metrizable, but we shall not use this fact.

\bigskip

\begin{prop}[{\bf Continuity}]
\label{P:continuity}
Let $G\acts (X,\mu)$ be an ergodic p.m.p actions of a \lcsc group $G$,
and $L$ be a group of class $\biinv$. 
Then the action of the centralizer $\Aut(X,\mu)^{G}$ on the space $Z^{1}(G\times X,L)$
of measurable cocycles $G\times X\to L$ is continuous.
\end{prop}
\begin{proof}[Proof for the discrete case]
We assume that both $G$ and $L$ are discrete groups, and denote them by 
$\Gamma$, $\Lambda$ respectively.
We should prove that given a finite set $F\subset \Gamma$ and an $\epsilon>0$ there is a 
a neighborhood $U$ of the identity in $\Aut(X,\mu)$ so that 
\begin{equation}
\label{e:cont-d-case}
	\mu\setdef{ x\in X}{\alpha(g,x)=\alpha(g,T(x))} >1-\epsilon\qquad (\forall g\in F)
\end{equation}
for all $T\in U \cap \Aut(X,\mu)^{\Gamma}$.
Indeed, for each $g\in \Gamma$ there exists a finite set $\Lambda_{g}\subset\Lambda$
so that 
\[
	\mu\setdef{x\in X}{\alpha(g,x)\in\Lambda_{g}}>1-\epsilon/2.
\]
Let $E_{g,\ell}=\setdef{x\in X}{\alpha(g,x)=\ell}$. 
Then there exists a neighborhood $U$ of the identity in $\Aut(X,\mu)$ so that
for $T\in U$ for each $g\in F$ we have:
\[
	\sum_{\ell\in\Lambda_{g}} \mu(E_{g,\ell}\symdiff T E_{g,\ell})<\epsilon/2.
\]
It now follows that (\ref{e:cont-d-case}) holds for all $T\in U$.
Restricting $T$ to commute with $\Gamma$ guarantees that $\alpha(g,T(x))$ is a cocycle.
\end{proof}
The general argument is only slightly more involved, and is postponed to the Appendix 
(section \ref{S:continuity}).
The following is a reformulation of Popa's \cite[Proposition 4.2]{Popa:cocycle}.
\begin{thm}[{\bf Local Rigidity of Cocycles}]
\label{T:local-rig}
Let $H< G$ be a closed subgroup which has relative property (T) in a \lcsc group $G$,
and let $G\acts (X,\mu)$ be a p.m.p. action such that the restricted action $H\acts (X,\mu)$ is ergodic.

Then for any group of finite type $L\in\fintyp$ the map 
\[
	Z^{1}(G\times X, L)\subset  Z^{1}(H\times X, L)\overto{} H^{1}(H\times X,L)
\]
has a discrete image. 
More precisely, there exists a compact subset $K\subset G$ and $\epsilon>0$ s.t.
for any two cocycles $\alpha,\beta :G\times X\to L$ with 
\[
	\mu\setdef{x\in X}{ d(\alpha(g,x),\beta(g,x))<\epsilon}>1-\epsilon\qquad (g\in K)
\]
there exists a measurable $\phi:X\to L$ so that as $H$-cocycles we have:
\begin{equation}
\label{e:abH}
	\alpha(h,x)=\phi(h.x)\,\beta(h,x)\,\phi(x)^{-1}\qquad(h\in H).
\end{equation}
\end{thm}
Here we shall the proof for the special case of \emph{discrete} target group, following 
Hjorth's \cite{Hjorth:converseDye:05}. 
This argument can be adapted to the case of $L$ being 
compact, or $L=\Inn(\Rel{R})$ for a ${\rm II}_{1}$-equivalence relation $\Rel{R}$. 
For the general case see Popa's paper \cite[Proposition 4.2]{Popa:cocycle}.
\begin{proof}[Proof - case of discrete target group]
Let $H<G$ and $G\acts (X,\mu)$ be as in the theorem, and let $\alpha:G\times X\to \Lambda$
be a measurable cocycle taking values in a discrete countable group $\Lambda$.
Using Theorem~\ref{T:alm-inv} we formulate the assumption that $H<G$ has 
relative property (T) in $G$ as follows:
there exists a compact subset $K\subset G$ and $\epsilon>0$ so that:
for any unitary $G$-representation $\pi:G\to U(\sH_{\pi})$ which has a unit vector
$v_{0}\in\sH_{\pi}$ with 
\[
	\inf_{g\in K}|\ip{\pi(g)v_{0}}{v_{0}}|\ge 1-\epsilon
\]
there exists an $H$-invariant unit vector $v\in\sH_{\pi}$ with $\|v-v_{0}\|<1/10$.

The claim is that any cocycle $\beta:G\times X\to \Lambda$ with  $(\alpha,\beta)\in V_{K,\epsilon}$,
i.e., such that
\[
	\mu\setdef{x\in X}{\alpha(g,x)=\beta(g,x)}>1-\epsilon\qquad (g\in K),
\]
relation (\ref{e:abH}) holds.

To see this, consider the product $\tilde{X}=X\times \Lambda$ with 
$\tilde{\mu}=\mu\times m_{\Lambda}$, 
where $m_{\Lambda}$ is the counting measure on $\Lambda$.
With $\alpha, \beta:G\times X\to\Lambda$ consider the measure preserving 
$G$-action on the infinite measure space $(\tilde{X},\tilde{\mu})$ by
\[
	g: (x,\ell)\mapsto (g.x,\, \alpha(g,x)\,\ell\,\beta(g,x)^{-1})
\]
and defines a unitary representation $\pi$ on $\sH=L^{2}(\tilde{X},\tilde{\mu})$.
Denoting by $F_0\in\sH$ the characteristic function of the set $X\times\{e\}\subset X\times \Lambda$, 
we note that $F_0$ is a unit vector with
\[
	\ip{\pi(g) F_0}{ F_0}=\mu\setdef{x\in X}{\alpha(g,x)=\beta(g,x)}>1-\epsilon\qquad(g\in K).
\]
Thus there exists an $H$-invariant unit vector $F\in\sH$ close to $ F_0$. 
For $x\in X$ consider $|F(x,\cdot)|^{2}$ as a function $\Lambda\to[0,\infty)$.
By Fubini for $\mu$-a.e. $x\in X$ this function is summable; let $w(x)$ denote its sum,
$p(x)$ denote its maximal value, and $k(x)$ -- the number of points on which this value is 
attained: 
\[
	p(x)=\max_{\ell\in L} |F(x,\ell)|^{2},\qquad
	k(x)={\rm card}\setdef{\ell\in  \Lambda}{|F(x,\ell)|^{2}=p(x)}.
\]
Then $w$, $p$ and $k$ are measurable $H$-invariant functions; hence by ergodicity 
these functions are $\mu$-a.e. constants: $w(x)=\|F\|^{2}=1$, 
$p(x)=p$ and $k(x)=k$.
In fact we claim that $k=1$. 
Indeed $k(x)\cdot p(x)\le w(x)$ gives $p\le 1/k$; in particular 
$|F(x,e)|^{2}\le 1/k$. Thus 
\[
	\| F_0-F\|^{2}\ge 1-1/\sqrt{k}
\] 
which would contradict $\|F-F_{0}\|<1/10$ unless $k=1$. 
Hence we can define a measurable function $\phi:X\to \Lambda$ by $|F(x,\phi(x))|^{2}=p$. 
As $F$ is $H$-invariant we get that for each $h\in H$:
\[
	\phi(h.x)=\alpha(h,x)\phi(x)\beta(h,x)^{-1}
\]
which is equivalent to (\ref{e:abH}).
\end{proof}
\begin{rem}
\label{R:general-rigidity}
In the general case of $L\in\fintyp$ (see \cite[Proposition 4.2]{Popa:cocycle}), where
$L$ is imbedded in the unitary group $U$ of a finite von Neumann algebra $M$, 
Popa applies the relative property (T) to the representation associated to $M$,
and then induced through $X$ to a unitary representation of $G$. 
After some work he shows that $\alpha:G\times X\to L<U$ 
can be untwisted to a homomorphism $\ro:G\to L<U$ 
(or reduced to a cocycle $\ro:G\times Y\to L<U$) by a conjugation in $U$.
It is here that Lemma~\ref{L:target} is used to untwist $\alpha$ to some $\ro':G\to L$
(or $\ro':G\times Y\to L$) within $L$ itself.
\end{rem}

\subsection{Proof of Popa's Cocycle Superrigidity Theorem}

Consider the fibered product $(Z,\zeta)=(X\times_{Y}X,\mu\times_{\nu}\mu)$
with the diagonal action of $G$. 
Let $C=\Aut(Z,\zeta)^{G}$ denote the centralizer of $G$ in $\Aut(Z,\zeta)$,
and let $C^{0}$ denote the connected component of the identity in $C$ 
considered with the weak topology.
By the relative malleability assumption, $C^{0}$ contains a map of the form 
\[
	F:(x,y)\mapsto (T(y),S(x))
\]
where $T,S\in \Aut(X,\mu)^{G}$.
The relative weak mixing assumption means that $H$ acts ergodically on $(Z,\zeta)$.

Given any measurable cocycle $\alpha:G\times X\to L$ where $L\in\fintyp$,
denote by $\alpha_{i}:G\times Z\to L$ the pull-back of $\alpha$ via the projections $p_{i}:Z\to X$, 
$p_{i}(x_{1},x_{2})=x_{i}$:
\[
	\alpha_{i}(g,(x_{1},x_{2}))\defq\alpha(g,x_{i})\qquad(i=1,2).
\]
The group $C^{0}< C=\Aut(Z,\zeta)^{G}$ acts continuously on the space 
$Z^{1}(G\times Z,L)$ of all measurable cocycles (Proposition~\ref{P:continuity}).
It follows from the Local Rigidity Theorem~\ref{T:local-rig}, that the space 
$H^{1}(H\times Z,L)$ inherits a discrete topology from $Z^{1}(G\times Z,L)$,
and therefore the action of the connected group $C^{0}$ on $Z^{1}(G\times Z,L)$
gives rise to a trivial action on $H^{1}(H\times Z,L)$.
This means  that any two cocycles in the $C^{0}$-orbit on $Z^{1}(G\times Z,L)$
are cohomologous as $H$-cocycles.
In particular, there exists a measurable map $F:Z=X\times_{Y}X\to L$ so that 
\[
	\alpha(h,x_{1})=F(h.x_{1},h.x_{2})\,\alpha(h,S(x_{2}))\, F(x_{1},x_{2})^{-1}\qquad(h\in H<G).
\]
By Theorem~\ref{T:untwisting} there exist a measurable map $\phi:X\to L$
and a measurable cocycle $\ro_{0}: H\times Y\to L$ so that:
\[
	\alpha(h,x)=\phi(h.x)\,\ro_{0}(h,p(x))\,\phi(x)^{-1}\qquad(h\in H).
\]
Define a cocycle  $\beta:G\times X\to L$  by 
$	\beta(g,x)\defq\phi(g.x)^{-1}\alpha(g,x)\phi(x) $
and let 
\[
	G_{1}\defq\setdef{g\in G}{\beta(g,x)\text{ is a.e. a function of  } y=p(x)}.
\]
It follows from the cocycle identity that $G_{1}$ is a subgroup of $G$,
and we already know that $H<G_{1}$. It remains to show that $G_{1}=G$.

If condition (b) in the theorem is satisfied, $H$ is w-normal in $G$.
Lemma~\ref{L:normal} states that for any subgroup $H<K<G_{1}$ we have $N_{G}(K)<G_{1}$
(the assumption that $H\acts (X,\mu)$ is weakly mixing relative to $(Y,\nu)$
passes to all the sup-groups $K>H$).
The assumption that $H$ is w-normal in $G$ implies, using transfinite induction
starting with $H_{0}=H$, that $G_{1}=G$.

If condition (b') in the theorem is satisfied, i.e., $H$ is only wq-normal, one
makes use of the second part of Lemma~\ref{L:normal}.
Let $H=H_{0}<H_{1}<\dots<H_{\eta}=G$ be the well ordered chain
of subgroups as in the definition of wq-normality.
For each ordinal $j\le \eta$ let $H'_{j}=\bigcup_{i<j} H_{i}$ 
and let $A_{j}$ be the set of elements $g\in G$ for which
the group $K_{j,g}\defq H_{j}'\cap g^{-1}H_{j}'g$ is not compact,
and hence acts on $(X,\mu)$ weakly mixing relative to $(Y,\nu)$.
By the assumption of wq-normality $H_{j}$ is generated by $A_{j}$.
Thus Lemma~\ref{L:normal} shows that if $H_{i}<G_{1}$ for all $i<j$, 
then $H_{j}<G_{1}$, and $G=G_{1}$ follows by the transfinite induction.

With $G=G_{1}$ shown, it follows that $\beta(g,x)=\ro(g,p(x))$ for some 
measurable function $\ro:G\times Y\to L$, which is necessarily a cocycle, 
and we conclude that 
\[
	\alpha(g,x)=\phi(g.x)\,\ro(g,p(x))\,\phi(x)^{-1}\qquad(g\in G)
\]
as claimed.
If $Y$ is a point space with the trivial action $\ro$ is just a homomorphism $G\to L$.
This completes the proof of Theorem~\ref{T:PopaCSR}.
\bigskip

\subsection{Malleable actions - proof of Theorem~\ref{T:malleable}}

\subsubsection*{Bernoulli Actions}
For reader's convenience we start with the generalized Bernoulli actions 
(see Popa \cite[Lemma 4.5]{Popa:cocycle}).
Let $I$ be a countable set with an action
$\sigma:G\to\Sym(I)$ of a \lcsc group, and let $(X,\mu)=(X_{0},\mu_{0})^{I}$ be the product
space with the $G$-action by coordinate permutation. 
The diagonal $G$-action on $(X\times X,\mu\times\mu)$ is
the Bernoulli action of $G$ on $(X_{0}\times X_{0},\mu_{0}\times \mu_{0})^{I}$.
Hence malleability of $G\acts (X,\mu)$ follows from the following observations
 \ref{L:Flip} and \ref{L:Delta}.
\begin{lem}
\label{L:Flip}
Let $(X_{0},\mu_{0})$ be a non-atomic probability measure space.
Then in the Polish group $\Aut(X_{0}\times X_{0},\mu_{0}\times\mu_{0})$,
considered with the weak topology, the flip $F_{0}(x,y)=(y,x)$ is in the (path)
connected component of the identity. 
\end{lem}
\begin{proof}
Without loss of generality we may assume that $(X_{0},\mu_{0})$ is the unit interval
$[0,1]$ with the Lebesgue measure $m$.
For $0\le t\le 1$ let $T_{t}\in\Aut([0,1]^{2},m^{2})$ be defined by
\[
	T_{t}(x,y)\defq \left\{ \begin{array}{ll} (x,y) &\text{if } x, y\in [0,t]\\
	(y,x) & \text{otherwise} \end{array} \right..
\]
Then $T_{0}=F_{0}$ is the flip, and $T_{1}=\Id$ is the identity.
The path 
\[
	t\in[0,1]\mapsto T_{t}\in\Aut([0,1]^{2},m^{2})
\] 
is continuous in the strong, and hence also in the weak, topologies.
\end{proof}
Let $(Z_{0},\zeta_{0})$ be some probability space, and $(Z,\zeta)=(Z_{0},\zeta_{0})^{I}$
be the product space with the Bernoulli action of $G$.
Let $\Delta:\Aut(Z_{0},\zeta_{0})\overto{}\Aut(Z,\zeta)$ denote the diagonal embedding.
The following fact is obvious:
\begin{lem}
\label{L:Delta}
$\Delta$ is a continuous embedding, and its image commutes with $G$.
\end{lem}
Since Bernoulli actions with a non-atomic base space are included in the family of Gaussian
actions (\ref{R:Bernoulli-is-Gaussian}), the above also follows from the following case.

\bigskip

We now turn to Gaussian actions (described in \ref{SS:Gaussian}) corresponding to
an orthogonal representation $\pi:G\to O(\sH)$.
First assume that $G$ is countable and the action is constructed from some $\pi(G)$-invariant
set $V=\{v_{i}\}_{i\in I}\subset \sH$; hence $G\acts (\RR^{I},\mu_{C})$.

The product measure $\mu_{C}\times\mu_{C}$ on $\RR^{I}\times\RR^{I}=(\RR^{2})^{I}$ 
is ``built from'' the standard Gaussian measures $\mu_{\RR^{2}}$ on $\RR^{2}$, and is
therefore invariant under the diagonal action of the rotation group $\SO(\RR^{2})$.
More precisely, the claim is that the rotation $R_{\theta}$ by angle $\theta$ acting on 
$\RR^{I}\times\RR^{I}$ by:
\[
	R_{\theta}:((x_{i})_{i\in I},(y_{j})_{j\in I})\mapsto ((x_{i}\cos{\theta} -y_{i}\sin{\theta})_{i\in I}, 
	(x_{j}\sin{\theta} +y_{j}\cos{\theta})_{j\in I})
\] 
preserves the product measure $\mu_{C}\times\mu_{C}$. 
Indeed the latter is the Gaussian measure $\mu_{\tilde{C}}$, where $\tilde{C}$ is a 
``block diagonal matrix'' $C\otimes I_{2}$ corresponding to the set
$\tilde V=V\otimes e_{1}\ \cup\ V\otimes e_{2}\subset \sH\otimes\RR^{2}\cong \sH\oplus\sH$.
For any finite subset $J\subset I$ the measure $\mu_{\tilde{C}_{J}}$ on $(\RR^{2})^{J}$
is invariant under the rotations $R_{\theta}$. 
Therefore $\mu_{\tilde{C}}=\mu_{C}\times\mu_{C}$ is also $R_{\theta}$-invariant.
The continuous path $\setdef{R_{\theta}}{0\le\theta\le\pi/2}$ in $\Aut(\RR^{I}\times\RR^{I},\mu_{C}\times\mu_{C})$
commuting with the diagonal $G$-action connects the identity with the map
\[
	R_{\pi/2}:(x,y)\mapsto(-y,x)
\]
proving the malleability of the Gaussian $G$-action on $(\RR^{I},\mu_{C})$.
This argument extends to Gaussian actions of non-discrete groups using a dense countable
subgroup $\Lambda<G$ and using the continuity argument as outlined in \ref{SS:Gaussian}.
This completes the proof of Theorem~\ref{T:malleable}.
\begin{rem}
This argument is in fact similar to a rotation argument used by Popa in the non-commutative
context in \cite[1.6.3]{Popa:RigidityI}.
\end{rem}

\section{Applications}
\label{S:Applications}

\begin{proof}[Proof of Theorem~\ref{T:superOE}]

First we reduce the weak morphism $X\supset X'\overto{\theta}Y'\subset Y$
to a morphism of the form  $X\overto{\theta'}Y'\subset Y$ (see \cite[Section 2]{Furman:OE:99} 
for a more detailed discussion).
First note that there exists a measurable ``retraction'' map 
$\pi:X\to X'$ with $(x,\pi(x))\in \Rel{R}_{X,\Gamma}$ for $\mu$-a.e. $x\in X$.
Indeed, enumerating elements of $\Gamma$ as $\{ g_{n}\}_{n\in\NN}$ 
(say, with $ g_{1}=e$) let $\pi(x)\defq g_{n(x)}.x$
where $n(x)=\min\setdef{ n\in\NN}{ g_{n}.x\in X'}$, 
which is finite for $\mu$-a.e. $x\in X$ by ergodicity.
Consider the map  
\[
	\sigma:X\overto{\pi} X'\overto{\theta} Y'\subset Y.
\]
Since $\pi_{*}\mu\sim \mu_{|X'}$ we have $\sigma_{*}\mu\sim\nu_{|Y'}\prec \nu$.

Observe that $\sigma(\Gamma.x)\subset\Lambda. \sigma(x)$ for $\mu$-a.e. $x\in X$,
and since $\Lambda$ is assumed to act freely on $(Y,\nu)$ (after possibly disregarding a null set),
this defines a measurable cocycle $\alpha:\Gamma\times X\to \Lambda$ by
\[
	\sigma( g.x)=\alpha( g,x).\sigma(x).
\]
Applying Popa's Cocycle Superrigidity Theorem~\ref{T:PopaCSR} there is a homomorphism
$\ro:\Gamma\to\Lambda$ and a measurable $\phi:X\to \Lambda$ so that
\[
       \alpha( g,x)=\phi( g.x)\,\ro( g)\,\phi(x)^{-1}.
\]
Define a measurable map $T:X\to Y$ by $T(x)\defq \phi(x)^{-1}.\sigma(x)$.
More explicitly, we have a measurable countable partition $X=\bigcup X_{\lambda}$
into sets $X_{\ell}=\phi^{-1}(\{\ell\})$, so that for $x\in X_{\ell}$:
$T(x)=\ell. \sigma(x)$.
In particular, this shows that $T_{*}\mu\prec\nu$. 
Let 
\[
	f(y)=\frac{dT_{*}\mu}{d\nu}(y),\quad
	Y_{1}=\setdef{y\in Y}{f(y)>0}, \quad
	F=f\circ T.
\]
Note that $F:X\to\RR$ is a $\Gamma$-equivariant measurable function, because
\begin{eqnarray*}
	T( g.x)&=&\phi( g.x)^{-1}.\sigma( g.x)
	=\phi( g.x)^{-1}\alpha( g, x). \sigma(x)\\
	&=&	\ro( g)\phi(x)^{-1}.\sigma(x)=\ro( g). T(x).
\end{eqnarray*}
Hence, by ergodicity,  $F(x)$ is $\mu$-a.e. a finite positive constant 
$c=\nu(Y_{1})/\mu(X)=\nu(Y_{1})$,
and the morphism $T:X\to Y_{1}\subset Y$ satisfies
\[
	T( g.x)=\ro( g).T(x)\qquad( g\in\Gamma).
\]
Since $T( g.x)=T(x)$ for every $ g\in\Gamma_{0}=\Ker(\ro)$, the map 
$T$ factors through $(X_{1},\mu_{1})$ -- the space of $\Gamma_{0}$-ergodic components. 
The $\Gamma$-action on $(X_{1},\mu_{1})$ factors through the action of
$\Gamma_{1}=\ro(\Gamma)\cong\Gamma/\Gamma_{0}$, and the corresponding map
$T_{1}:(X_{1},\nu_{1})\to(Y_{1},\nu_{1})$ is $\Gamma_{1}$--equivariant.
Being a quotient of an ergodic action $\Gamma\acts(Y_{1},\nu_{1})$ is ergodic.
This proves (1)--(3).

\medskip

For the proof of (i)--(iii) we need the following
\begin{clm}
There exist positive measure subsets $A\subset X$ and $B\subset Y_{1}$ 
and fixed elements 
$a\in\Gamma$, $b\in\Lambda$, so that $T(A)=B$ and
$	T(x)=b.\theta(a.x)$ for a.e. $x\in A$.
\end{clm}
\begin{proof}
This follows from the construction of the map $T:X\to Y_{1}$ as
\[
	T(x)=\ell_{x}\cdot\theta(\pi(x))=\ell_{x}.\theta( g_{n(x)}.x)
	\qquad
	\text{where}\qquad \ell_{x}=\phi(x)^{-1}.
\]
Indeed, since both $\Gamma$ and $\Lambda$ are countable, there
exists a positive measure set $A\subset X$ so that $n(x)$ is constant $n_{0}$
on $A$ (we denote $a= g_{n_{0}}$), 
while $\ell_{x}$ is constant $b$ on $A$.
The set $B=b.\theta( a.A)=T(A)$ has positive measure
and $A=T^{-1}(B)$.
\end{proof}

(i). We assume the action $\Gamma\acts (X,\mu)$ to be free, and $\theta:X'\to Y'$
to be an injective morphism of the restrictions $\Rel{R}_{X,\Gamma}|_{X'\times X'}$ into
$\Rel{R}_{Y,\Lambda}|_{Y'\times Y'}$. 
Using the above claim it follows that also $T:A\to B$ is an \emph{injective} relation morphism 
of the restrictions $\Rel{R}_{X,\Gamma}|_{A\times A}\to\Rel{R}_{Y,\Lambda}|_{B\times B}$.
Note that for any $g\in\Gamma$, if $x\in A\cap g^{-1} A$ then  
\[
	\Rel{R}_{X,\Gamma}|_{A\times A}\ni (x, g. x)\ \overto{T}\   
	(T(x),\ro( g).T(x))\in \Rel{R}_{Y,\Lambda}|_{B\times B}.
\]
Thus the fact that $T:A\to B$ is an injective relation morphism, together
with the freeness of the $\Gamma$-action, imply that for 
$ g\in\Ker(\ro)\setminus\{e\}$ we have $\mu(A\cap  g^{-1} A)=0$.
Therefore
\[
	|\Gamma_{0}|\cdot \mu(A)=\sum_{ g\in\Gamma_{0}}\mu( g A)
	=\mu(\bigcup_{  g\in\Gamma_{0}} g A)\le\mu(X)=1
\]
implying the finiteness of  $\Gamma_{0}\normal\Gamma$.

On the other hand, if $\Gamma_{0}$ is finite, then $(X_{1},\mu_{1})$ is just the
space of $\Gamma_{0}$-orbits in $(X,\mu)$, and there exists a Borel section
$f:X_{1}\to X$ of the finite extension $(X,\mu)\to (X_{1},\mu_{1})=(X,\mu)/\Gamma_{0}$.   
Observe that $T$ mapping $X'\defq f(X_{1})$ to $Y_{1}$ as
\[
	X\supset X'\overto{f^{-1}} X_{1}\overto{T_{1}}Y_{1}\subset Y 
\]
is an injective relation morphism from 
$\Rel{R}_{X,\Gamma}|_{X'\times X'}$ into 
$\Rel{R}_{Y,\Lambda}|_{Y_{1}\times Y_{1}}= \Rel{R}_{Y_{1},\Gamma_{1}}$.

\medskip
 
(ii). We assume that $\nu(Y)<\infty$ and that $\theta:X'\to Y'$
is a surjective morphism of the restrictions  
$\Rel{R}_{X,\Gamma}|_{X'\times X'}$ onto
$\Rel{R}_{Y,\Lambda}|_{Y'\times Y'}$.
Using the above claim it follows that $T:A\to B$ is an surjective relation morphism 
of the restrictions $\Rel{R}_{X,\Gamma}|_{A\times A}$ onto 
$\Rel{R}_{Y,\Lambda}|_{B\times B}$.

Let $\ell\in\Lambda$ be such that $\nu(B\cap\ell^{-1}B)>0$.
Then for $\nu$-a.e. $y\in B\cap\ell^{-1}B$ we have 
$(y,\ell.y)\in \Rel{R}_{Y,\Lambda}|_{B\times B}$ 
and, by surjectivity, we can find and $x\in A$ and $g\in \Gamma$ 
with $y=T(x)$ and $\ell.y=T(g.x)=\ro(g).T(x)$.
The essential freeness of the action $\Lambda\acts (Y,\nu)$ implies
that $\ell=\ro(g)\in\Gamma_{1}$.

Reversing this implication, we conclude that $\nu(B\cap \ell^{-1}B)=0$ whenever 
$\ell\not\in\Gamma_{1}$.
Let $\{\ell_{i}\}_{i\in I}$ be representatives of distinct cosets 
$\Lambda/\Gamma_{1}$. 
Thus $\ell_{i}^{-1}\ell_{j}\not\in\Gamma_{1}$ for $i\neq j$ and 
$\nu(\ell_{i}B\cap\ell_{j}B)=\nu(B\cap\ell_{i}^{-1}\ell_{j}B)=0$. Hence
\[
	[\Lambda:\Gamma_{1}]\cdot \nu(B)=\sum_{i\in I}\nu(\ell_{i}B)
	=\nu(\bigcup_{i\in I}\ell_{i}B)\le\nu(Y)<\infty
\]
which yields $[\Lambda:\Gamma_{1}]<\infty$.

Let $\Lambda'=\bigcap_{\ell\in\Lambda} \ell\Gamma_{1}\ell^{-1}$  -- the maximal
normal subgroup in $\Lambda$ contained in $\Gamma_{1}$.
It has finite index in $\Lambda$. 
Thus $\Gamma'=\ro^{-1}(\Lambda')$ has finite index in $\Gamma$,
and acts ergodically on $(X,\mu)$, by the aperiodicity assumption on
$\Gamma\acts (X,\mu)$.
Hence $\Lambda'\acts (Y_{1},\nu_{1})$, which is a quotient of $\Gamma'\acts (X,\mu)$,
is an ergodic action.
Since $\Lambda'\normal\Lambda$ the decomposition of $(Y,\nu)$ into 
$\Lambda'$-ergodic components (one of which is $(Y_{1},\nu_{1})$)
is acted upon by the finite group $\Lambda/\Lambda'$.
This action is transitively due to ergodicity of the $\Lambda$-action on $(Y,\nu)$,
and each has the same $\nu$-measure.
Thus $\nu(Y)/\nu(Y_{1})=k$ is an integer dividing the index $[\Lambda:\Gamma_{1}]$.
 
\medskip
 
(iii) follows from the combination of (i) and (ii). The last claim, that if $c(\theta)=1$ then 
the weak relation morphism $\theta$ extends to a relation morphism bijection 
$\tilde\theta:X\to Y$ (up to null sets), follows from an easy adaptation  of
\cite[Proposition 2.7]{Furman:OE:99}. The latter treats weak relation isomorphisms rather
than weak bijective morphisms.

\end{proof}

\begin{proof}[Sketch of the proof of Theorem~\ref{T:ME+}]
The case of a single cocycle $A=\{\alpha\}$ is basically \cite[Theorems B,C]{Furman:OE:99}.
For a set $A$ of cocycles we take $Y_{A}=\bigvee_{\alpha\in A} Y_{\{a\}}$.
In other words $L^{\infty}(Y_{A},\nu_{A})$ is the abelian sub-algebra of $L^{\infty}(X,\mu)$
generated by  $L^{\infty}(Y_{\{a\}},\nu_{\{a\}})$, $\alpha\in A$.
The properties (1) and (2) of $(Y_{A},\nu_{A})$ are evident.
The structure of $(Y_{A},\nu_{A})$ when 
$A=\{\alpha_{1},\dots,\alpha_{n}\}$ is a finite set
follows from the fact (see \cite{Furman:Outer:05}) that 
Ratner's theorem allows to explicitly describe 
$\Gamma$-invariant ergodic probability 
measures on $\prod_{i=1}^{n}Y_{\{\alpha_{i}\}}=\prod_{i=1}^{n} \bar{G}/\Lambda_{i}$
with marginals $m_{\bar{G}/\Lambda_{i}}$ as homogeneous measures for 
some subdiagonal embedding of $\bar{G}^{k}$ in $\bar{G}^{n}$.
\end{proof}

\begin{proof}[Proof of Theorem~\ref{T:alg-not-quotient}]
Starting from a weak relation morphism $\theta$ from $\Rel{R}_{X,\Gamma}$
to $\Rel{R}_{H/\Delta,\Lambda}$, Theorem~\ref{T:superOE} gives a 
homomorphism $\ro_{1}:\Gamma\to\Lambda$ and the chain of $\Gamma$-equivariant 
measurable quotients 
\[
	T:(X,\mu)\overto{erg}(X_{1},\mu_{1})\overto{T_{1}}(Y_{1},\nu_{1}),
	\qquad Y_{1}\subset H/\Delta.
\]
The image $\Gamma_{1}=\ro_{1}(\Lambda)$ acts ergodically on a non-trivial
space $(Y_{1},\nu_{1})$. In fact $\nu_{1}\prec m_{H/\Delta}$ is absolutely continuous 
with respect to Haar measure.
Thus $\Gamma_{1}<\Delta<H$ is infinite and hence unbounded, and therefore acts
ergodically on $(H/\Delta,m_{H/\Delta})$ by Howe-Moore theorem.
So $Y_{1}=Y=H/\Delta$ and $\nu_{1}=\nu=m_{H/\Delta}$.
Let  $Z_{1}=Z_{H}(\Gamma_{1})$ denote the centralizer of $\Gamma_{1}$ in $H$.
For later use we observe that $Z_{1}$ is a proper algebraic subgroup in $H$,
and therefore is a null set with respect to the Haar measure $m_{H}$.

Let $f:H/\Delta\to H$ be a Borel cross-section of the projection $H\to H/\Delta$,
and let $c:H\times H/\Delta\to\Delta$ be the corresponding measurable cocycle 
\[
	c(h,y)=f(h.y)\,h\,f(y)^{-1}\qquad(h\in H,\ y\in H/\Delta).
\]
Denote by $\alpha:\Gamma\times X\to\Delta$ the lift of $c$ to $X$, namely
$\alpha( g,x)\defq c(\ro_{1}( g),T(x))$. We have
\[
	\alpha( g,x)=\Phi_{1}( g.x)\,\ro_{1}( g)\,\Phi_{1}(x)^{-1},
	\qquad\text{where}\qquad
	\Phi_{1}:X\overto{T} H/\Delta\overto{f} H.
\]
The assumption that $\Gamma\acts(X,\mu)$ is $\Delta$-Cocycle Superrigid 
means that there exists a homomorphism $\ro_{2}:\Gamma\to\Delta$ and a 
measurable map $\Phi_{2}: X\to \Delta$ so that
\[
	\alpha( g,x)=\Phi_{2}( g.x)\,\ro_{2}( g)\,\Phi_{2}(x)^{-1}.
\]
We shall deduce a contradiction by comparing $\Phi_{1}$ to $\Phi_{2}$.
To do so rewrite the identity 
\[
	\Phi_{2}( g.x)\,\ro_{2}( g)\,\Phi_{2}(x)^{-1}
	=\alpha( g,x)=
	\Phi_{1}( g.x)\,\ro_{1}( g)\,\Phi_{1}(x)^{-1}
\] 
as $\Phi( g.x)=\ro_{1}( g)\,\Phi(x)\ro_{2}( g)^{-1}$ with $\Phi(x)\defq \Phi_{1}(x)^{-1}\Phi_{2}(x)$.
Next define a measurable map $\Psi:X\times X\to H$ by
$\Psi(x',x)\defq \Phi(x')\Phi(x)^{-1}$. 
We have
\[
	\Psi( g.x', g.x)=\ro_{1}( g)\,\Psi(x',x)\,\ro_{1}( g)^{-1}.
\]
Therefore the push-forward probability measure $\omega=\Psi_{*}(\mu\times\mu)$ on $H$
is invariant under conjugations by $\Gamma_{1}=\ro_{1}(\Gamma)<\Delta<H$.

It is a general fact (see \cite[p. 38]{Furman:MM:01}) that the only probability measure 
on a semi-simple center free group $H$ which is invariant under conjugations 
by elements of a subgroup $\Gamma_{1}<H$ is supported on the centralizer 
$Z_{1}=Z_{H}(\Gamma_{1})$.

Thus $\Phi(x')\Phi(x)^{-1}\in Z_{1}$ for $\mu\times\mu$-a.e. $(x',x)\in X\times X$.
By Fubini, there exists $h_{0}\in H$ (take $h_{0}=\Phi(x')^{-1}$ for $\mu$-a.e. $x'$)
so that for $\mu$-a.e. $x\in X$: 
\[
	\Phi_{2}(x)^{-1}\Phi_{1}(x)=\Phi(x)^{-1}\in h_{0}Z_{1}.
\] 
Since $\Phi_{2}:X\to\Delta$ it follows that for $\mu$-a.e. $x\in X$:
\[
	\Phi_{1}(x)\in \Phi_{2}(x)h_{0}Z_{1}\subset 
	\bigcup_{h\in\Delta} h h_{0}Z_{1}.
\]
The right hand side is a countable union of translates of $Z_{1}$, and therefore
is a null set with respect to the Haar measure $m_{H}$ on $H$.
At the same time ${\Phi_{1}}_{*}\mu\prec m_{H}$ because
\[
	{\Phi_{1}}_{*}\mu=f_{*}m_{H/\Delta}
\] 
is the restriction of the Haar measure on $H$ to a fundamental domain for $\Delta<H$.
The contradiction shows that there could not exist a weak relation morphism from
$\Rel{R}_{X,\Gamma}$ to $\Rel{R}_{H/\Delta,\Lambda}$.
\end{proof}


\begin{proof}[On Remark~\ref{R:Quotients-of-Bernoulli}(2)]
Let $\Gamma$ be a higher rank lattice, $(X,\mu)=(X_{0},\mu_{0})^{I}$ an ergodic
generalized Bernoulli action, and $(X,\mu)\overto{}(Y,\nu)$ a measurable
quotient where each $g\in\Gamma$ acts with finite entropy.
Then $(Y,\nu)$ is a point. This follows from the following facts: 
\begin{enumerate}
\item 
By Stuck and Zimmer \cite{StuckZimmer:94} an ergodic action 
of a higher rank lattice is essentially free, unless the action is 
transitive one on a finite set.
Since $\Gamma\acts (X,\mu)$ is weakly mixing, it follows that
$\Gamma\acts (Y,\nu)$ is essentially free. 
\item 
Higher rank lattices $\Gamma$ always contain a copy of $\ZZ^{2}$; 
Hence restricting to $\ZZ^{2}$ we get that $\ZZ^{2}\acts (Y,\nu)$ is an essentially
free action which is a quotient of a Bernoulli action of $\ZZ^{2}$ on $(X,\mu)$.
\item
Hence $\ZZ^{2}\acts (Y,\nu)$ is isomorphic to a non-trivial Bernoulli action 
$\ZZ^{2}\acts (Z_{0},\zeta_{0})^{\ZZ^{2}}$, being an essentially free quotient of
a Bernoulli action $\ZZ^{2}\acts (X,\mu)$.
\item
The generators $g_{1},g_{2}$ of the free Abelian group $\ZZ^{2}$,
must have infinite infinite Kolmogorov-Sinai entropy on 
$(Y,\nu)=(Y_{0},\nu_{0})^{\ZZ^{2}}$, because these are Bernoulli $\ZZ$-actions with 
a non-atomic base space 
$(Y_{1},\nu_{1})=(Y_{0},\nu_{0})^{\ZZ}$:
\[
	h(Y,\nu,g_{1})=H(Y_{1},\nu_{1})=\infty\cdot H(Y_{0},\nu_{0})=\infty
\]
since $(Y_{0},\nu_{0})$ is not a singleton.
\end{enumerate}
\end{proof}

\begin{proof}[Proof of Proposition~\ref{P:extensions}]
If $(X',\mu')\overto{p}(X,\mu)$ is not relatively weakly mixing, then there exists
an intermediate extension $p:(X',\mu')\overto{r}(Z,\zeta)\overto{q}(X,\mu)$,
where $q$ is a non-trivial relatively compact extension, i.e.,
$(Z,\zeta)=(X\times K/K_{0},\mu\times m_{K/K_{0}})$
and the $G$-action on $Z$ is given by a measurable cocycle $\alpha:G\times X\to K$.
Using $\cpt$-Cocycle Superrigidity, $\alpha(g,x)=\phi(g.x)\ro(g)\phi(x)^{-1}$
with $\ro:G\to K$ homomorphism, and $\phi:X\to K$ is a measurable map.
Then $f:Z\to Z$ given by $f(x,kK_{0})=(x,\phi(x)kK_{0})$ conjugates
the given action $G\acts (Z,\zeta)$ to the diagonal action $G\acts X\times K/K_{0}$,
$g:(x,kK_{0})\mapsto (g.x,\ro(g)kK_{0})$.
The map 
\[
	X'\overto{r}Z\overto{f}X\times K/K_{0}\overto{p_{2}}K/K_{0}
\]
shows that $G\acts (X',\mu')$ has the isometric action $G\acts K/K_{0}$ as a quotient,
contrary to the assumption that $G\acts (X',\mu')$ is weakly mixing. 
\end{proof}

\begin{proof}[Proof of Theorem~\ref{T:semi-direct-rel}]
Define a measurable cocycle $\Gamma\times X\to\Inn(\Rel{S})$
as follows: 
for  $g\in\Gamma$ and $\mu$-a.e. $x\in X$ let $\alpha_{g,x}\in\Inn(\Rel{S})$ 
be defined by
\[
	((x,y),(g.x,\alpha_{g,x}.y))\in\Rel{Q}.
\]
It follows that $\alpha$ is a cocycle. Since $\Inn(\Rel{S})\in\fintyp$ and
$\Gamma\acts (X,\mu)$ is assumed to be $\fintyp$-Cocycle Superrigid we have
a homomorphism $\ro:\Gamma\to\Inn(\Rel{S})$ and a measurable map
$\phi:X\to\Inn(\Rel{S})$ so that 
$\alpha(g,x)=\phi(g.x)\,\ro(g)\,\phi(x)^{-1}$.
Let $f\in\Aut(X\times Y,\mu\times\nu)$ be given by $f(x,y)=(x,\phi(x).y)$.
This is an inner relation automorphism of $\Rel{R}\times\Rel{S}$.

Note that the diagonal $\Gamma$-action on $(X\times Y,\mu\times\nu)$, 
$g:(x,y)\mapsto (g.x, \ro(g.y)$ is also inner
to $\Rel{R}\times\Rel{S}$, and $\Phi$ identifies the orbit relation of this $\Gamma$-action 
with $\Rel{Q}$.
\end{proof}

\section{Appendix}
\label{S:continuity}

\subsection{Proof of Proposition~\ref{P:continuity} - the general case}
We consider a fixed measurable cocycle $\alpha:G\times X\to L$ where $L$
is a Polish group with some bi-invariant metric $d$. 
For $T\in\Aut(X,\mu)^{G}$ we shall measure the distance between $\alpha$ and 
the shifted cocycle $\alpha(g,T(x))$ by the function:
\[
	f_{T}(g,x)= d\left(\alpha(g,x), \alpha(g,T(x))\right)
\]
and will use the following sets depending on the parameters $\epsilon>0$, $\delta>0$:
\begin{eqnarray*}
	E_{\epsilon}(T,g)&=&\setdef{x\in X}{f_{T}(g,x)>\epsilon},\\
	A_{\epsilon,\delta}(T)&=&\setdef{g\in G}{ \mu(E_{\epsilon}(T,g))<\delta}.
\end{eqnarray*}
\begin{lem}
The above functions and sets satisfy the following relations: 
\begin{eqnarray}
	f_{T}(gh,x)&\le& f_{T}(g,h.x)+f_{T}(g,x),\qquad(g,h\in G)	\label{e:subadd}\\
	f_{TS}(g,x)&\le& f_{T}(g,S(x))+f_{S}(g,x)\qquad(T,S\in \Aut(X,\mu)^{G})	\label{e:TS}\\
	\mu(E_{\epsilon_{1}+\epsilon_{2}}(T,gh))&\le & 
		\mu(E_{\epsilon_{1}}(T,h))+\mu(E_{\epsilon_{2}}(T,g))
				\label{e:Egh}\\
	\mu(E_{\epsilon_{1}+\epsilon_{2}}(TS,g))&\le & 
		\mu(E_{\epsilon_{1}}(T,g))+\mu(E_{\epsilon_{2}}(S,g))
				\label{e:ETS}\\
	A_{\epsilon_{1},\delta_{1}}(T)^{-1}\cdot A_{\epsilon_{2},\delta_{2}}(T)
		&\subset& A_{\epsilon_{1}+\epsilon_{2},\delta_{1}+\delta_{2}}(T)
				\label{e:AA}
\end{eqnarray}
For  fixed $g\in G$, $\epsilon>0$ and $\delta>0$ there exist 
an open neighborhood $U=U(g,\epsilon,\delta)$
of the identity in $\Aut(X,\mu)$ so that:
\begin{equation}
\label{e:ETU}
	g\in \bigcap_{T\in U\cap \Aut(X,\mu)^{G} }A_{\epsilon,\delta}(T).
\end{equation}
\end{lem}
\begin{proof}
(\ref{e:subadd}). 
This inequality follows from the cocycle identity, combined with the following 
``parallelogram inequality'' which takes advantage of the bi-invariance of the 
metric $d$ on $L$:
\begin{eqnarray*}
	d(a_{1} a_{2},\, b_{1} b_{2})&=&d(b_{1}^{-1}a_{1},\, b_{2}a_{2}^{-1})
	\le d(b_{1}^{-1}a_{1},e)+d(e, b_{2}a_{2}^{-1})\\
	&=&d(a_{1},b_{1})+d(a_{2}, b_{2})
\end{eqnarray*}
for every $a_{1}, a_{2}, b_{1}, b_{2}\in L$. 

(\ref{e:TS}) is immediate from the triangle inequality
\[
	d(\alpha(g,x),\alpha(g,TSx))\le d(\alpha(g,x),\alpha(g,Sx))
		+d(\alpha(g,Sx),\alpha(g,TSx)).
\] 

(\ref{e:Egh}) follows from the inclusion 
\[
	E_{\epsilon_{1}+\epsilon_{2}}(T,gh)\subset  
		E_{\epsilon_{1}}(T,h)+ h^{-1}E_{\epsilon_{2}}(T,g)
\]
which is implied by (\ref{e:subadd}).

(\ref{e:ETS}) is immediate from (\ref{e:TS}).

(\ref{e:AA}) follows from (\ref{e:Egh}) and the easy fact:
$A_{\epsilon,\delta}(T)=A_{\epsilon,\delta}(T)^{-1}$.

(\ref{e:ETU}). 
Choose an $\epsilon/2$-dense sequence $\{\ell_{i}\}_{i=1}^{\infty}$ in $L$, let 
\[
	X_{i}=\setdef{x\in X}{ d(\alpha(g_{0},x),\ell_{i})<\epsilon/2 },
\]
and choose $N$ large enough to ensure $\mu(\bigcup_{i=1}^{N}X_{i})>1-\delta/2$.
Set
\[
	U(g,\epsilon,\delta)=\bigcap_{i=1}^{N}\setdef{T\in\Aut(X,\mu)}
	{\mu(TX_{i}\symdiff X_{i})<\frac{\delta}{2N}}.
\]
For $T\in U(g,\epsilon,\delta)$ we have 
$\mu(\bigcup_{i=1}^{N} (X_{i}\cap T^{-1}X_{i}))>1-\delta$,
while for $x\in X_{i}\cap T^{-1}X_{i}$:
\[
	d(\alpha(g,x), \alpha(g,T(x)))\le
	d(\alpha(g,x),\ell_{i})+d(\ell_{i}, \alpha(gT(x)))
	<\frac{\epsilon}{2}+\frac{\epsilon}{2}=\epsilon.
\]
This completes the proof of the Lemma.
\end{proof}

\begin{proof}[Proof of Proposition~\ref{P:continuity}]
Fix a compact subset $K\subset G$ and $\epsilon>0$ and $\delta>0$.

It follows from (\ref{e:ETU}) that there exists an open neighborhood $U_{0}$
of the identity in $\Aut(X,\mu)$ so that there is a set $B=B_{\epsilon/3,\delta/3}$ 
of positive Haar measure of elements $g\in G$ with 
\[
	B\subset \bigcap_{T\in U_{0}\cap\Aut(X,\mu)^{G}} A_{\epsilon/3,\delta/3}(T).
\]
Then, by a standard Fubini argument, there exits an open neighborhood $V$ of 
the identity in th e\lcsc group $G$ s.t. $V\subset B^{-1}B$. 
Hence for every $T\in U_{0}\cap\Aut(X,\mu)^{G}$ we have, using (\ref{e:AA}):
\[
	V\subset B^{-1}B\subset A_{\epsilon/3,\delta/3}(T)^{-1}A_{\epsilon/3,\delta/3}(T)
	\subset A_{2\epsilon/3,2\delta/3}(T).
\]
Since $K$ is compact, there exists a finite set $\{g_{1},\dots,g_{m}\}\subset K$ so that 
$K\subset \bigcup_{i=1}^{m} g_{i}V$.
For each $g_{i}$ there is a neighborhood $U_{i}$ of the identity in $\Aut(X,\mu)$
so that $g_{i}\in A_{\epsilon/3,\delta/3}(T)$ for all $T\in U_{i}\cap \Aut(X,\mu)^{G}$.
Therefore taking $U=U_{0}\cap U_{1}\cap\dots\cap U_{m}$, and using (\ref{e:AA})
again, we can deduce that 
\[
	K\subset \bigcap_{T\in U\cap \Aut(X,\mu)^{G}} A_{\epsilon,\delta}(T).
\] 
This proves the Proposition.
\end{proof}

\nocite{Popa:cocycle} \nocite{Popa:RigidityI} \nocite{Popa:RigidityII}
\nocite{Fernos:relT}
\nocite{Furstenberg:Szemeredi:77}
\nocite{Furman:ME:99} \nocite{Furman:OE:99} \nocite{Furman:Outer:05}
\nocite{Hjorth:converseDye:05}
\nocite{HjorthKechris:MAMS:05}
\nocite{Vaes:afterPopa:06}
\nocite{Zuk:GAFA:03}
\nocite{Shalom:MGT:05}
\nocite{OrnsteinWeiss:BernoulliAmen:87}



\def\cprime{$'$}
\providecommand{\bysame}{\leavevmode\hbox to3em{\hrulefill}\thinspace}
\providecommand{\href}[2]{#2}

\end{document}